\documentclass[review,1p,times]{elsarticle}

\journal{arXiv}

\biboptions{authoryear}

\usepackage{amsmath}
\usepackage{amssymb}
\usepackage{graphicx}
\usepackage{adjustbox}
\usepackage[colorlinks=true]{hyperref}
\usepackage{subcaption}
\usepackage{tikz}

\usetikzlibrary{
    decorations.pathmorphing, 
    arrows, 
    positioning,
    shapes.geometric
}
\tikzset{
    ->,
    >=stealth',
    shorten >= 1pt,
    auto,
    node distance = 1cm, thick,
    bellman/.style={circle, draw},
    square/.style={regular polygon sides = 4, draw},
    textnode/.style={}
}
\newcommand{\drawsquiggle}[2]{
    \draw[
        -stealth,
        decoration = {
            snake,
            amplitude = 0.4mm,
            segment length = 1.5mm,
            post length = 0.9mm},
        decorate
    ] (#1) + (#2) -- (#1);
}
\newcommand{\drawHDsquiggle}[2][0.5]{\drawsquiggle{#2}{-#1, #1};}

\begin{document}

\begin{frontmatter}

\title{Capacity planning of renewable energy systems using stochastic dual dynamic programming}

\author[epoc,ntnu,nve]{J. Hole\corref{cor1}}
\ead{jarand.hole@ntnu.no}
\author[epoc]{A.B. Philpott}
\author[epoc]{O. Dowson}
\address[epoc]{Electric Power Optimization Centre, University of Auckland, New Zealand}
\address[ntnu]{Norwegian University of Science and Technology, Norway}
\address[nve]{Norwegian Water Resources and Energy Directorate, Norway}
\cortext[cor1]{Corresponding author}

\begin{abstract}

We present a capacity expansion model for deciding the new electricity generation and transmission capacity to complement an existing hydroelectric reservoir system. The objective is to meet a forecast demand at least expected cost, namely the capital cost of the investment plus the expected discounted operating cost of the system. The optimal operating policy for any level of capacity investment can be computed using stochastic dual dynamic programming. We show how to combine a multistage stochastic operational model of the hydro system with a capacity expansion model to create a single model that can be solved by existing open-source solvers for multistage stochastic programs without the need for customized decomposition algorithms. We illustrate our method by applying it to a model of the New Zealand electricity system and comparing the solutions obtained with those found in a previous study.
\end{abstract}

\begin{keyword}
OR in energy; stochastic dual dynamic programming; capacity expansion; multistage stochastic programming
\end{keyword}

\end{frontmatter}

\newpage

\section{Introduction}

Electricity systems around the world are transitioning to technologies with zero or near-zero carbon emissions \citep{ipccwgiii}. The most popular sources of renewable energy (wind and solar) are either intermittent or diurnal and so some form of flexibility is required to match supply and demand when the wind is not blowing or the sun is not shining. One option is thermal peaker plants, which can quickly ramp up and down to balance supply and demand. Another option is to store and discharge energy from some form of storage, which may be medium-term storage in the form of hydro-reservoirs, or short-term storage in the form of batteries. Since future wind and solar generation are uncertain, the optimal operation of storage is the solution to a stochastic control problem of some complexity. A third option is to invest in interconnectors, through which neighbouring power systems can share excess electricity and flexibility. 

In addition to the intermittency challenges, isolated hydro-dominated energy systems like Iceland, Brazil, and New Zealand, face a ``dry-year'' risk, in which a sustained period of low inflows leaves insufficient water in the reservoirs to meet demand. To compensate, the system must either be willing to periodically shed demand at high cost (and social interruption), build and operate thermal peakers as a back-up source of energy, or over-build intermittent renewable capacity to ensure sufficient supply. In each case, extra transmission may be needed to move power from the newly built generation to demand loads. The main contribution of this paper is a model which can be used to correctly size the investment and analyse the operation of each option.

A key challenge when analysing the operation of a hydro-dominated system is the trade-off between the value of using the water stored in a reservoir in the current period, compared with the value of keeping the water for a future period. Most commonly, the operation of the system is modelled as a multistage stochastic program and solved using stochastic dual dynamic programming (SDDP) \citep{pereira1991}. One such example is the Brazilian energy system, which has used SDDP to plan its operation for over 25 years \citep{maceira2008ten,maceiral2018twenty}. Because of the use of both hydro and thermal generation, the problem is often called the \textit{hydro-thermal scheduling problem}.

A number of authors have augmented multistage stochastic operational models with capacity expansion decisions. As a brief sampling of the literature, we point to the following works \citep{newham2008,rebennack2014,wu2016two,optgen,bruno2016,thome2019stochastic,lara2020electric,bodal2022,genx}. The literature can be categorized into two broad streams, which we refer to as \textit{static} and \textit{dynamic} investment models.

In \textit{static} investment models, the capacity investment and operation are treated as a two-stage problem, in which the agent chooses a level of new capacity investment, and then evaluates the operation of the energy system with the investment complete. Static models focus on the ideal end-state of the energy system, but not the timing or sequence of decisions that are required to get there. In its simplest form, a static model enumerates all potential investments and evaluates their capital and operating costs before choosing the best. More sophisticated static models treat the operation of the energy system as the subproblem of a two-stage model that has first-stage capacity investment decisions (modelled by binary or continuous variables) optimized using a form of decomposition. The second-stage is solved, often using SDDP, to guide the solution of the first stage. A unifying feature of static models is that they decompose the investment and operational decisions, and that their operational problems model a finite time horizon. The use of a finite-horizon model to evaluate the steady-state operational cost of the system forces models to use excessively long horizons to mitigate end-of-horizon effects. Moreover, because these models use decomposition, the authors each code a customized algorithm to solve their particular problem.

\textit{Dynamic} investment models co-optimize the capacity investment and operation of the energy system over time. In contrast to static models, changes in capacity can occur more than once over the time horizon. Dynamic models are usually formulated as monolithic multistage stochastic programs, typically with integer decisions for the investments. Examples of such models are \cite{dominguez2016investing}, \cite{liu2017multistage}, and \cite{backe2022empire}. Because of the difficulty in solving large instances, various approximations or heuristics are used to generate sub-optimal solutions. A unifying feature of dynamic models is that they model a finite time horizon.

In this paper, we propose a static capacity investment model that co-optimizes the strategic level investment of new capacity to minimize the steady-state operational cost of how the energy system should be operated in the face of uncertain inflows and quantity of variable renewable energy. Our model is based on the policy graph modelling framework of \cite{dowson_policy_graph} and can be solved using existing open-source multistage stochastic programming solvers such as SDDP.jl \citep{dowson2021sddp} without the need for the practitioner to code their own specialized decomposition algorithm. The single, co-optimized,  infinite-horizon model formulation is the main novel contribution of this paper. We apply our model to a case study that investigates the capacity expansion of variable renewable energy in the New Zealand electricity system. We consider different investment scenarios, and to verify that the method works, we compare our optimal solutions with those computed by enumeration by \citet{philpott2023onslow}. The ability to solve our model without the need to implement custom decomposition algorithms is particularly important, because it means that our approach can easily be applied by other researchers to other energy systems and case studies. Thus, this paper also serves as a demonstration of the flexibility of the policy graph modelling approach and of SDDP.jl.

The paper is laid out as follows. In the next section we describe the proposed method for the integration of capacity investments in SDDP.jl. In Section~\ref{sec:wind} we describe how we model wind investment in SDDP.jl. In Section~\ref{sec:case} we give a short description of the New Zealand case study, and compare the results from our investment model and previous work in Section~\ref{sec:results}, before we conclude in Section~\ref{sec:conclusions}. \ref{appendix:jade} gives the full formulation of the SDDP stage problem and \ref{appendix:system} provides parameters for the generation plants modelled in the case study. More detailed results from the case study are presented in \ref{appendix:storage} and \ref{appendix:watervalue}.

\section{Modelling investments using a policy graph}\label{sec:theory}

In this section we introduce the hydro-thermal scheduling problem, explain how to model the finite- and  infinite-horizon variants as a policy graph, and expand the policy graph to include investment decisions.

The classical version of the hydro-thermal scheduling model is a finite-horizon discrete-time stochastic optimal control problem. The {\em stages} of this problem are indexed $1,2, \ldots, T$. The {\em states} of this problem are reservoir levels measured at the end of each stage (denoted $X$) and {\em actions} (denoted $U$) are water releases through generating turbines, flows through spill and river arcs, dispatchable non-hydro generation, and load reduction. In each stage, each  reservoir $r$ experiences a random inflow $\omega_r (t)$. The upper-case notation used for $X$ and $U$ indicate that these are also random variables.

The general formulation of this hydro-thermal scheduling problem (HTP) is as follows:
\begin{equation*}
\begin{array}{r l l}
\textbf{HTP}: \min& {\mathbb E}[\sum_{t=1}^{T} C_t(X(t),U(t))] \\
\text{s.t.} & X(t)= f_t(X(t-1),U(t),\omega(t))  &t=1,2, \ldots, T\\
& X(0) = \bar{x} \\
& U(t) \in \mathcal{U}_t(x(t-1), \omega(t))     & t=1,2, \ldots, T \\
& X(t) \in \mathcal{X}_t                        & t=1,2, \ldots, T.
\end{array}
\end{equation*}

Here, $ C_t(X(t),U(t))$ is the cost of meeting demand (possibly with load shedding) in stage $t$ from dispatchable generation. The finite horizon ignores actions after stage $T$, but future costs and constraints on $X(T)$ can be modelled by suitable choices of $C_T(X(T),U(T))$ and $\mathcal{X}_T$. Note that these exogeneous modelling choices will affect the optimal solution of \textbf{HTP}. We revisit this issue below.

The transition function $f_t(X(t-1),U(t),\omega(t))$ maps the incoming state $X(t-1)$ at the beginning of stage $t$ to a random outgoing state $X(t)$ at the end of stage $t$. For example:
$$X(t) = f_t(X(t-1),U(t),\omega(t)) = X(t-1) + B \cdot U(t) +\omega(t),$$
where the matrix $B$ maps the control flows $U(t)$ to the change in reservoir level. The vector of inflows $\omega(t)$ is assumed to have a finite probability distribution (drawn from the sample space $\Omega(t)$) that is independent of that in any previous stage.  Observe that random variables $X$ and $U$ are required to be measurable with respect to the history of the random inflow process, so they must satisfy {\em non-anticipative} constraints that we have suppressed in the formulation.

\begin{figure}[!ht]
    \centering
        \begin{tikzpicture}
        \node[bellman]   (a)  []              {};
        \node[square,node distance=2cm] (b)  [right of = a] {$\mathbf{SP}_{t = 1}$};
        \node[textnode,node distance=2cm] (c)  [right of = b] {$\mathbf{...}$};
        \node[square,node distance=2cm] (d)  [right of = c] {$\mathbf{SP}_{t = T}$};
        \path[out=0,in=180]
        	(a)	 edge [] node [above] {} (b)
            (b)	 edge [] node [above] {} (c)
            (c)	 edge [] node [above] {} (d);
        \drawHDsquiggle[0.65]{b}
        \drawHDsquiggle[0.65]{d}
        \end{tikzpicture}
    \caption{The policy graph structure for \textbf{HTP}.}
    \label{fig:policy-graph1}
\end{figure}
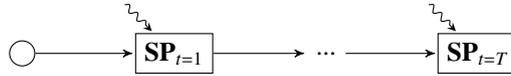

The model \textbf{HTP} can be defined by a \textit{policy graph} as described in \citet{dowson_policy_graph}. The policy graph defines the dynamic structure of the decision problem we are modelling, showing how actions affect the state variables and how realizations of the random variables are revealed over time.
The policy graph that represents \textbf{HTP} is shown in Figure~\ref{fig:policy-graph1}. 

There are $T$ stages with subproblems $\mathbf{SP}_t$, where each $t$ is a time step when decisions are made. The random variable $\omega(t)$ is shown as squiggly lines in the policy graph. Each subproblem $\mathbf{SP}_t$ has the following form:

\begin{equation*}
\begin{array}{r l}
\mathbf{SP}_t(x(t-1), \omega(t)): \min\limits_{u(t), x(t)} & C_t( x(t), u(t), \omega(t)) \\
\text{s.t.} &   x(t) = f_t(x(t-1),u(t),\omega(t)) \\
 & u(t) \in  \mathcal{U}_t(x(t-1), \omega(t)) \\
 & x(t) \in \mathcal{X}_t.
 \end{array}
\end{equation*}

The need to choose a suitable terminal cost function $C_T$ can be mitigated by solving the discounted infinite-horizon problem:
\begin{equation*}
\begin{array}{r l l}
\mathbf{HTP-\infty}: \min& {\mathbb E}[\sum_{t=1}^{\infty} \rho^t C_t(X(t),U(t))] \\
\text{s.t.} & X(t)= f_t(X(t-1),U(t),\omega(t))  & t=1,2, \ldots\\
& X(0) = \bar{x} \\
& U(t) \in \mathcal{U}_t(x(t-1), \omega(t))     & t=1,2, \ldots\\
& X(t) \in \mathcal{X}_t                        & t=1,2, \ldots.
\end{array}
\end{equation*}

\begin{figure}[ht]
    \centering
        \begin{tikzpicture}
        \node[bellman]   (a)  []              {};
        \node[square,node distance=2cm] (b) [right of = a] {$\mathbf{SP}_{t = 1}$};
        \node[textnode,node distance=2cm] (c) [right of = b] {$\mathbf{...}$};
        \node[square,node distance=2cm] (d)  [right of = c] {$\mathbf{SP}_{t = T}$};
        \path[out=0,in=180]
            (a)	 edge [] node [below] {} (b)
            (b)	 edge [] node [below] {$\rho$} (c)
            (c)	 edge [] node [below] {$\rho$} (d);
            \draw (d) to [out=-90,in=-90,looseness=0.5] node [below] {$\rho$}(b);
        \drawHDsquiggle[0.65]{b}
        \drawHDsquiggle[0.65]{d}
        \end{tikzpicture}
    \caption{The policy graph structure for $\mathbf{HTP-\infty}$.}
    \label{fig:stepwise-policy-graph}
\end{figure}
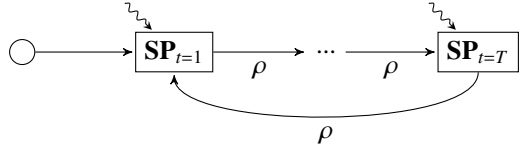

 Infinite-horizon discounted cost models can be represented by policy graphs that contain cycles. An example is shown in Figure~\ref{fig:stepwise-policy-graph}. In each step, the probability of transitioning from node $t$ to $t+1$ (or from node $T$ to node $1$) is $\rho$ where $\rho < 1$. Most commonly, $T$ is chosen so that one loop of the graph represents one year, and the nodes within each loop allow for seasonality in demand and inflows over the course of each year. Importantly, the subproblems $\mathbf{SP}_t$ do not change from the finite-horizon model, only the structure of the graph changes.

We can extend the $\mathbf{HTP-\infty}$ problem to include an investment variable $u_{inv}$ as follows:
\begin{equation*}
\begin{array}{r l l}
\mathbf{INV-HTP-\infty}: \min& {c_{inv}}^\top u_{inv} + {\mathbb E}[\sum_{t=1}^{\infty} \rho^t C_t(X(t),U(t))] \\
\text{s.t.} & X(t)= f_t(X(t-1),U(t),\omega(t))  & t=1,2, \ldots\\
& X(0) = \bar{x} \\
& U(t) \in \mathcal{U}_t(x(t-1), u_{inv}, \omega(t))     & t=1,2, \ldots\\
& X(t) \in \mathcal{X}_t(u_{inv})                     & t=1,2, \ldots.
\end{array}
\end{equation*}
Here $u_{inv}$ is the vector of new capacity that is installed before the start of stage $t=1$. The newly installed capacity modifies the feasibility set of the state variables $\mathcal{X}_t$, and it may also modify the feasibility set of the control variable $\mathcal{U}_t$. Here $c_{inv}$ is the overnight investment cost $I$ per unit of $u_{inv}$ [\$/MW], assuming that the invested capacity has an infinite lifetime. If the lifetime is finite, say $\tau$ years, we assume that the capacity is reinvested every $\tau$ years and set:
\begin{equation*}
    c_{inv} = I + \sum\limits_{i=1}^{\infty} I \cdot \beta^{i \cdot \tau} = \frac{I}{1-\beta^{\tau}}.
\end{equation*}

\begin{figure}[ht]
    \centering
        \begin{tikzpicture}
        \node[bellman]   (a)  []              {};
        \node[square,node distance=2cm] (b)  [right of = a] {$\mathbf{SP}_{inv}$};
        \node[square,node distance=2cm] (c)  [right of = b] {$\mathbf{SP}_{t = 1}$};
        \node[textnode,node distance=2cm] (d)  [right of = c] {$\mathbf{...}$};
        \node[square,node distance=2cm] (e)  [right of = d] {$\mathbf{SP}_{t = N}$};
        \path[out=0,in=180]
        	(a)	 edge [] node [below] {} (b)
            (b)	 edge [] node [below] {} (c)
            (c)	 edge [] node [below] {$\rho$} (d)
            (d)	 edge [] node [below] {$\rho$} (e);
            \draw (e) to [out=-90,in=-90,looseness=0.5] node [below] {$\rho$}(c);
        \drawHDsquiggle[0.65]{c}
        \drawHDsquiggle[0.65]{e}
        \end{tikzpicture}
    \caption{The policy graph structure for $\mathbf{INV-HTP-\infty}$.}
    \label{fig:investment-policy-graph}
\end{figure}
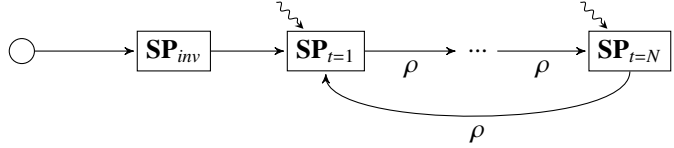

The policy graph in Figure~\ref{fig:stepwise-policy-graph} can be amended to accommodate capacity expansion decisions at the beginning of the infinite operating horizon. This gives the policy graph shown in Figure~\ref{fig:investment-policy-graph}. The decisions in $\mathbf{SP}_{inv}$ determine the system capacities, and the decisions in $\mathbf{SP}_t$ are operating decisions optimizing the control of the resulting system over an infinite horizon. 
 
In order to pass the capacity decisions from $\mathbf{SP}_{inv}$ into the cyclic policy graph, we add a state variable, $x_{inv}$, which is set by the investment control variable $u_{inv}$. Subproblem $\mathbf{SP}_{inv}$ is as follows:
\begin{equation*}
\begin{array}{r l}
\mathbf{SP}_{inv}: \min\limits_{u_{inv}}& {c_{inv}}^\top u_{inv} \\
\text{s.t.} & u_{inv} \ge 0 \\
  & x_{inv}(0) = u_{inv}
\end{array}
\end{equation*}

The new $x_{inv}$ state variable changes the constraints of each subproblem $\mathbf{SP}_t$. For dispatchable plants, investment increases the maximum output in any dispatch period. The same holds for grid investments, where the maximum power flow on a transmission line is increased. The effect of investment on the output of intermittent capacity such as wind or solar power is less straightforward.  In the next section, we describe a methodology that approximates wind generation from invested wind capacity for use in a medium term hydro-scheduling model.

\section{Modelling renewable investments in the subproblems}\label{sec:wind}

In the previous section we showed how capacity investments can be modelled
in SDDP.jl using a policy graph. In this section we show how to translate
this investment decision into capacity constraints on generation in each
stage problem. When generation comes from intermittent sources the
constraints are stochastic. We illustrate how we approximate these
constraints using wind as an example. The same method will apply for solar power.

In each stage problem, wind generation $w_{m}(h,t)$ (MW from wind farm $m$
in hour $h$ for week $t$ of the year) is a decision variable that
contributes to the electricity supply with a short-run marginal cost equal
to $0$. The output of wind farm $m$ in week $t$ and hour $h$ depends upon
its capacity $x_{m}^{inv}$, but also the (random) amount of wind available
at that time. This can be represented using by a random factor $\lambda
_{m,h,t}(\omega _{t})$ that limits the output through:
$$w_{m}(h,t,\omega _{t})\leq \lambda _{m,h,t}(\omega _{t})\cdot x_{m}^{inv},$$
where the inequality constraint gives that curtailment of wind generation is
possible at zero cost.

In our model we approximate the distribution of $\lambda _{m,h,t}(\omega _{t})$ by an empirical distribution from $Y$ years
of historical wind generation data. The wind farms for which these data are
available are grouped into geographical regions $r$, and $\lambda
_{m,h,t}(\omega_{t})$ is assumed to be the same value $\lambda
_{r,h,t}(\omega_{t})$ for each wind farm $m$ in region $r$. In any region $r$ and hour $h$ in week $t$, $\lambda_{r,h,t}(\omega_{t})$ can take on one of $Y$ values with equal
probability. Each value in this distribution (corresponding to year $y\in
\{1,2,\ldots ,Y\}$ ) is obtained by dividing the observed wind generation in
year $y$ in hour $h$ in week $t$ in region $r$ by the installed wind
capacity in that region. This model assumes wind generation in any hour is
independent of that in the previous hour. This ignores the short-term
temporal dependence of wind generation from hour to hour that might become
important if investment in short-term storage were being considered. The
standard approach (see e.g. \cite{liu2017multistage}) in this case is to use ``representative
days'' of wind generation, which model the dynamics at the
expense of introducing some perfect foresight into the decision making.

Since each weekly stage problem in the SDDP model is to be solved many
times, hourly demand and hourly wind generation are approximated by
piecewise constant load duration curves. The standard approach to estimating
these assumes a level of wind capacity $K_{r}$ for each region $r\in 
\mathcal{R}$ and subtracts the wind generation $\lambda _{r,h,t}(\omega
_{t})K_{r}$ from the demand $D_{r}(h,t)$ in that region to give a net demand:
$$\hat{D}_{r}(h,t,K_{r},\omega _{t})=D_{r}(h,t)-\lambda _{r,h,t}(\omega
_{t})K_{r},$$
for region $r$ in hour $h$ in week $t$ corresponding to the wind in historical year $%
\omega _{t}$. Here $\lambda _{r,h,t}(\omega _{t})$ is the load factor
computed from historical wind generation in region $r$ divided by capacity
in year $\omega _{t}$. Then total system net demand $\sum_{r\in \mathcal{R}}%
\hat{D}_{r}(h,t,K_{r},\omega _{t})$ is sorted to give a duration curve that
is approximated by load blocks $b=1,2,\ldots ,B$, each having constant
net load equal to average net load over the hours and historical years in
its block.

It is easy to see that different choices of wind capacity $K_{r}$ will
produce different values of net demand. Moreover the duration curve obtained
by sorting $\sum_{r\in \mathcal{R}}\hat{D}_{r}(h,t,K_{r},\omega _{t})$ will
sort the observations from different hours and historical years into a
possibly different order for each choice of $K_{r}$. A brute-force approach
to capacity optimization would enumerate all potential choices of $K_{r}$, $%
r\in \mathcal{R}$, and solve an SDDP model with this choice. In our\ model,
the values of $K_{r}$, $r\in \mathcal{R}$ are represented by state variables 
$x_{inv}$, which means that they are optimized by SDDP rather than
enumerated ex-ante. This requires a different approach.

The first approximation that we make is to assume that $K_{r}\,$\ is a fixed
proportion $\alpha _{r}$ of system wind capacity $K$, a single state
variable that will be optimized by SDDP. This implies that net demand in
region $r$ is:
$$\hat{D}_{r}(h,t,\alpha _{r}K,\omega _{t})=D_{r}(h,t)-\lambda _{r,h,t}(\omega
_{t})\alpha _{r}K,$$
and total net demand for wind capacity $K$ is:
$$\sum_{r\in \mathcal{R}}\hat{D}_{r}(h,t,\alpha _{r}K,\omega _{t})=\sum_{r\in 
\mathcal{R}}D_{r}(h,t)-\sum_{r\in \mathcal{R}}\alpha _{r}\lambda
_{r,h,t}(\omega _{t})K.$$

Given a nominal level $\bar{K}$ of wind capacity we compute $\sum_{r\in 
\mathcal{R}}\hat{D}_{r}(h,t,\alpha_r\bar{K},\omega _{t})$ and sort to give a
duration curve that is approximated by load blocks $b=1,2,\ldots ,B$.
Each load block $b$ for week $t$ contains $T(t,b)$ observations. We fix this
number for each load block.

To model the load duration curve as $K$ varies from $\bar{K}$, we make a
second approximation. To understand how $K$ affects the duration curve it is
helpful to fix week $t$ and write: 
\begin{eqnarray*}
\hat{D}(h,t,K,\omega _{t}) &=&\sum_{r\in \mathcal{R}}\hat{D}_{r}(h,t,\alpha
_{r}\bar{K},\omega _{t})-\sum_{r\in \mathcal{R}}\alpha _{r}\lambda
_{r,h,t}(\omega _{t})(K-\bar{K}) \\
&=&X(h,\omega _{t})-W(h,\omega _{t})(K-\bar{K}),
\end{eqnarray*}%
where $X(h,\omega _{t})=\sum_{r\in \mathcal{R}}\hat{D}_{r}(h,t,\alpha _{r}%
\bar{K},\omega _{t})$ and $W(h,\omega _{t})=\sum_{r\in \mathcal{R}}\alpha
_{r}\lambda _{r,h,t}(\omega _{t})$. Consider first the highest load block $%
b=1$ for week $t$, containing the $T(1,t)$ largest values of $X(h,\omega
_{t})$. As $K$ increases from $\bar{K}$, the indices $(h,\omega _{t})\in $ $b
$ will remain the same but the average load will decrease until a swap
occurs when we reach a critical value $K_{1}$, such that for $K>K_{1}$ and
some index $(h^{i},\omega _{t}^{i})\in b$:
$$X(h^{i},\omega _{t}^{i})-W(h^{i},\omega _{t}^{i})(K-\bar{K})<X(h^{j},\omega
_{t}^{j})-W(h^{j},\omega _{t}^{j})(K-\bar{K})$$
for $(h^{j},\omega _{t}^{j})\notin b$, and we swap $(h^{j},\omega _{t}^{j})$
into block 1 and $(h^{i},\omega _{t}^{i})$ out of block 1. When this happens
we will have:
$$X(h^{i},\omega _{t}^{i})-W(h^{i},\omega _{t}^{i})(K_{1}-\bar{K}%
)=X(h^{j},\omega _{t}^{j})-W(h^{j},\omega _{t}^{j})(K_{1}-\bar{K})$$
so $\hat{D}(h,t,K,\omega _{t})$ is continuous in $K$. (Observe that $X$ and $%
W$ are not continuous in $K$, since they can both jump at $K_{1}$.) Since
wind investment reduces net demand, we have $\hat{D}(h^{i},t,K,\omega
_{t}^{i})$ is decreasing in $K$, and:
$$X(h^{i},\omega _{t}^{i})-W(h^{i},\omega _{t}^{i})(K_{1}-\bar{K}%
)=X(h^{j},\omega _{t}^{j})-W(h^{j},\omega _{t}^{j})(K_{1}-\bar{K})$$
implies:
$$W(h^{i},\omega _{t}^{i})>W(h^{j},\omega _{t}^{j}).$$

This implies that $\sum_{(h,\omega _{t})\in 1}\hat{D}(h,t,K,\omega _{t})$ is
a convex decreasing piecewise linear function of $K$. This enables us to
model the effect of wind in block $b=1$ using linear inequality constraints
(which is required for SDDP). Unfortunately the above analysis does not
apply to the other blocks. Indeed it can be shown that for the lowest demand
block, $B(t)$, we have $\sum_{(h,\omega _{t})\in B(t)}\hat{D}(h,t,K,\omega
_{t})$ is a \emph{concave} decreasing piecewise linear function of $K$,
which cannot be represented in SDDP using linear inequality constraints.

Our approximation for week $t$ takes the load blocks constructed for $%
\sum_{r\in \mathcal{R}}\hat{D}_{r}(h,t,K,\omega _{t})$ for five different
values of $K$ ranging from 1 to 5 GW, where each load block $b$ is
restricted to contain exactly $T(t,b)$ observations. (We set $Y=1$ in our
experiments meaning that $\omega _{t}$ is ignored and $%
\sum_{b=1}^{B}T(t,b)=168$.) The observations (hours) that are in block $b$
for capacity choice $K$ are denoted $b(K)$. The average extra wind
generation in this block is 
\begin{eqnarray*}
&&\frac{\sum_{h\in b(K)}W(h,\omega _{t})}{T(t,b)}K-\frac{\sum_{h\in
b}W(h,\omega _{t})}{T(t,b)}\bar{K} \\
&=&\sum_{r\in \mathcal{R}}\alpha _{r}\left( \frac{\sum_{h\in b(K)}\lambda
_{r,h,t}(\omega _{t})}{T(t,b)}K-\frac{\sum_{h\in b}\lambda _{r,h,t}(\omega
_{t})}{T(t,b)}\bar{K}\right) 
\end{eqnarray*}%
Recall that $\frac{\sum_{h\in b(K)}\lambda _{r,h,t}(\omega _{t})}{T(t,b)}K$
is not continuous in $K$, since it can jump up when $b(K)$ changes. We
construct a continuous linear approximation by fitting a straight line
through $(K,$ $\frac{\alpha _{r}K\sum_{h\in b(K)}\lambda _{r,h,t}(\omega
_{t})}{T(t,b)})$ for five different values of $K$ ranging from 1 to 5 GW.
The slope $\mu _{r,b,t}$ of this line gives an estimate of the change in
wind generation in region $r$ for a unit increase in national wind capacity.
We can can then model wind generation $w_{r}(t,b)$ in block $b$ and week $t$
in SDDP using the investment state $x_{inv}$ and the constraint:
$$w_{r}(t,b)\leq \mu _{r,b,t}x_{inv}.$$

\begin{figure}[htbp]
\centering
\begin{subfigure}[b]{.32\textwidth}
\includegraphics[width=\textwidth]{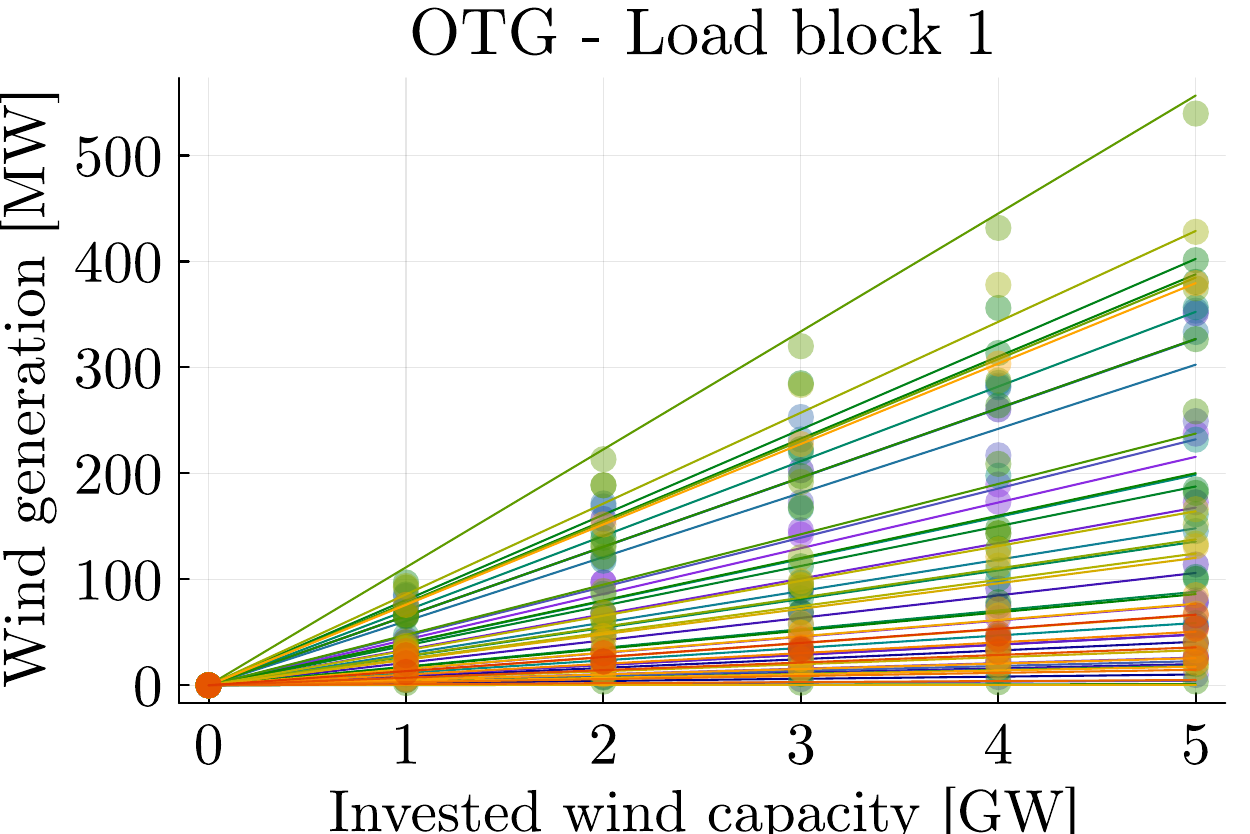}
\caption{Linear wind representation in Otago for load block 1}
\label{fig:otg1}
\end{subfigure} 
\hfill 
\begin{subfigure}[b]{.32\textwidth} 
\includegraphics[width=\textwidth]{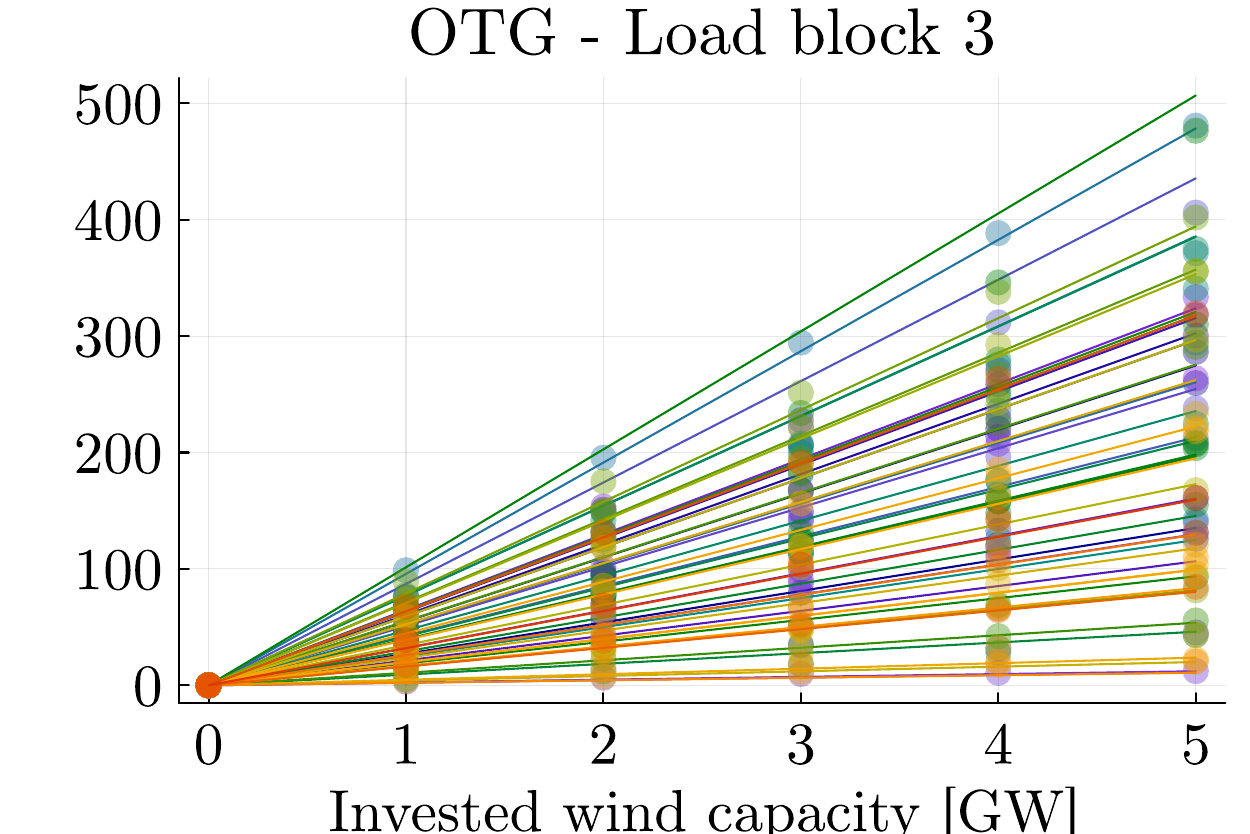} 
\caption{Linear wind representation in Otago for load block 3} 
\label{fig:otg3} 
\end{subfigure}
\hfill 
\begin{subfigure}[b]{.32\textwidth} 
\includegraphics[width=\textwidth]{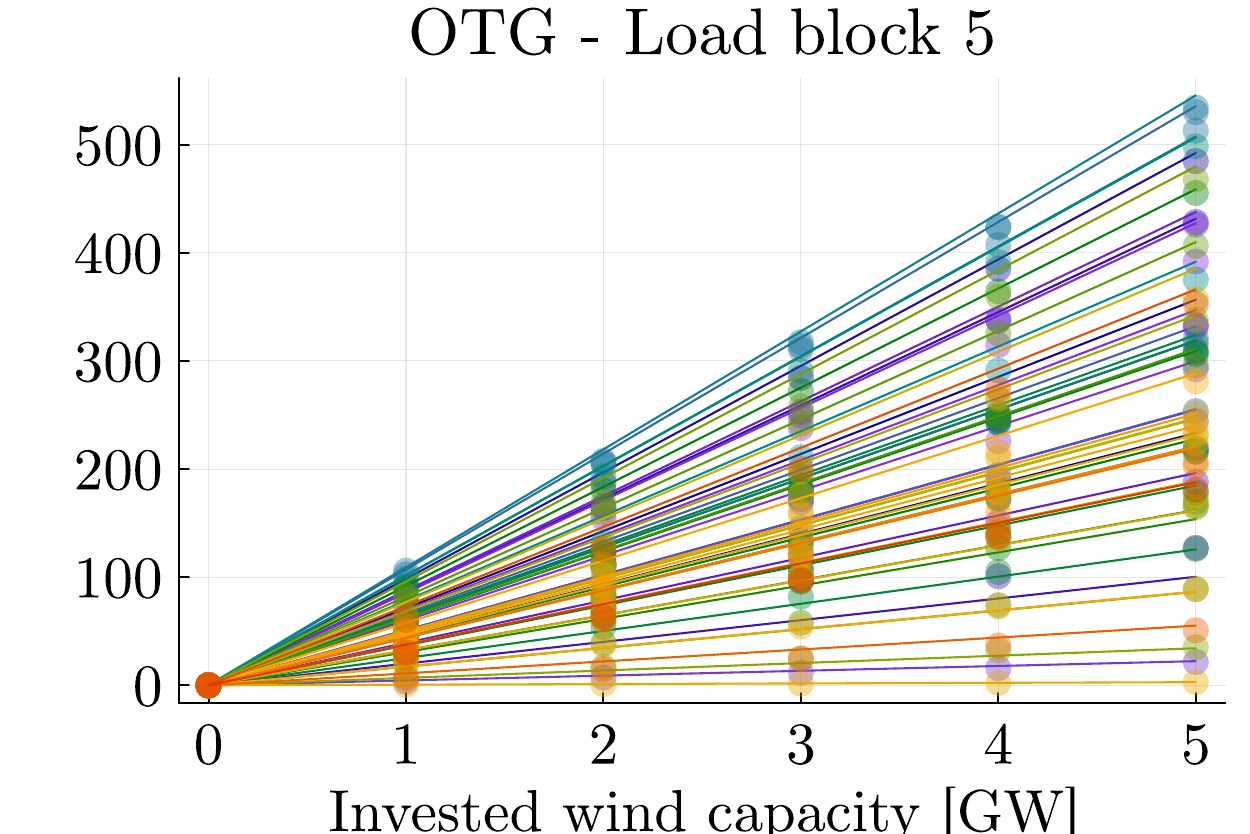} 
\caption{Linear wind representation in Otago for load block 5} 
\label{fig:otg5} 
\end{subfigure}
\caption{Otago wind generation (MW) versus $K$ (GW), where load block sizes $T(b,t)$ are determined using $\bar{K}$=2.5 GW. Different colours correspond to different weeks of the year.} 
\label{fig:windotg}
\end{figure}

\begin{figure}[!ht]
\centering
\begin{subfigure}[b]{.32\textwidth}
\includegraphics[width=\textwidth]{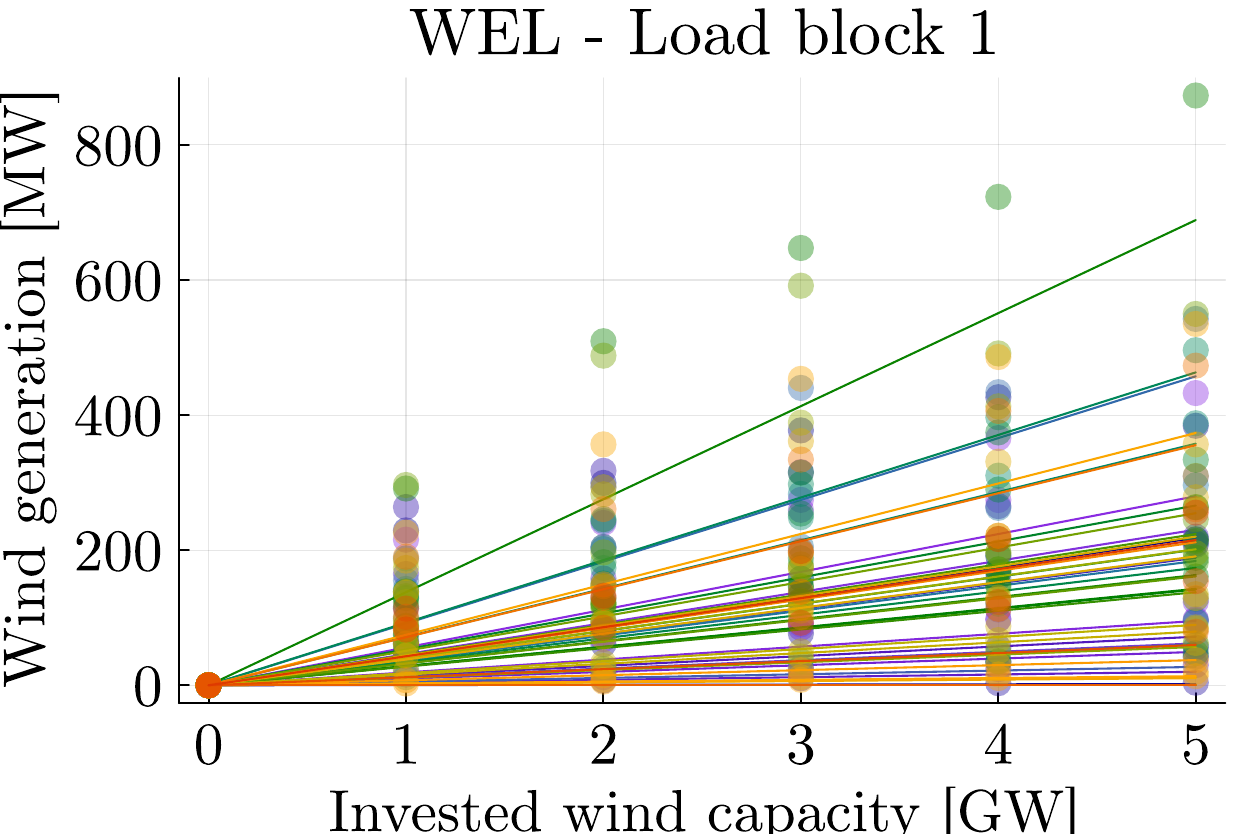}
\caption{Linear wind representation in Wellington for load block 1}
\end {subfigure} 
\hfill 
\begin{subfigure}[b]{.32\textwidth} 
\includegraphics[width=\textwidth]{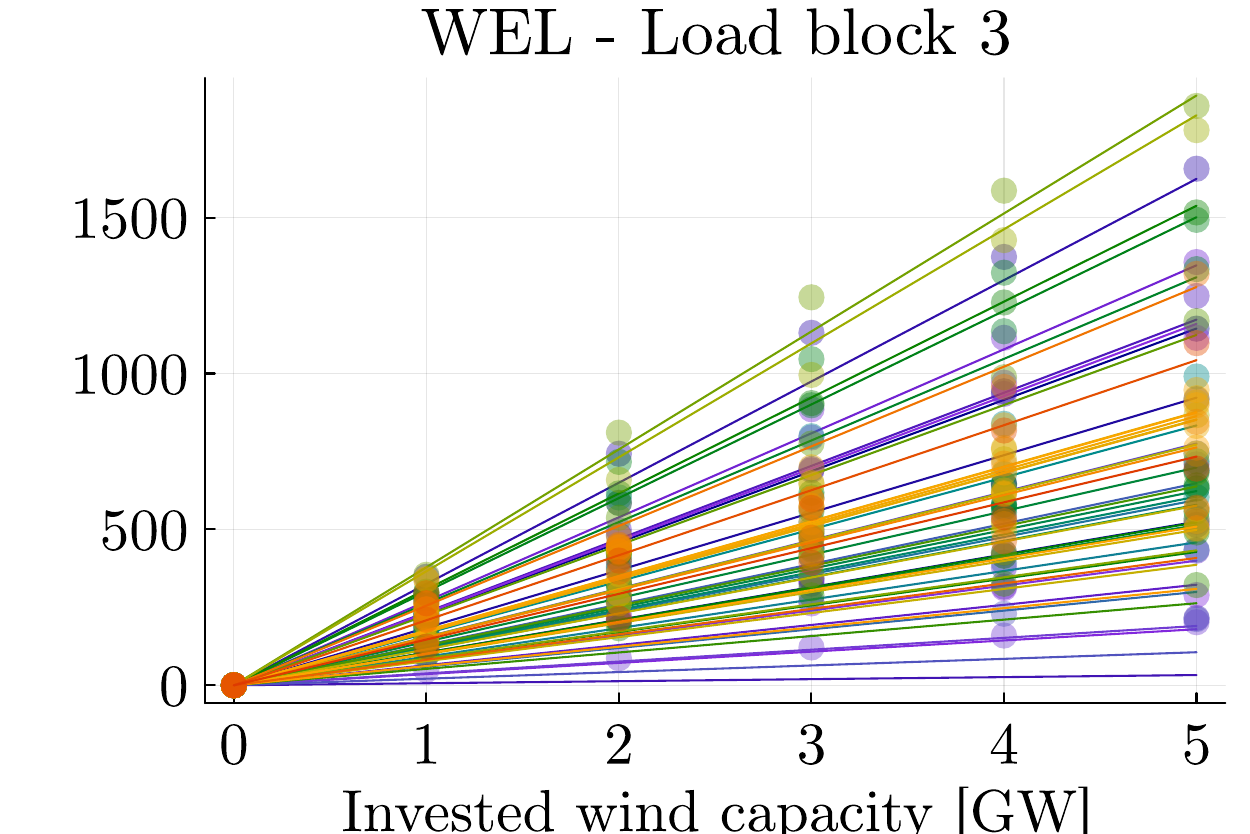}
\caption{Linear wind representation in Wellington for load block 3}
\end{subfigure}
\hfill 
\begin{subfigure}[b]{.32\textwidth} 
\includegraphics[width=\textwidth]{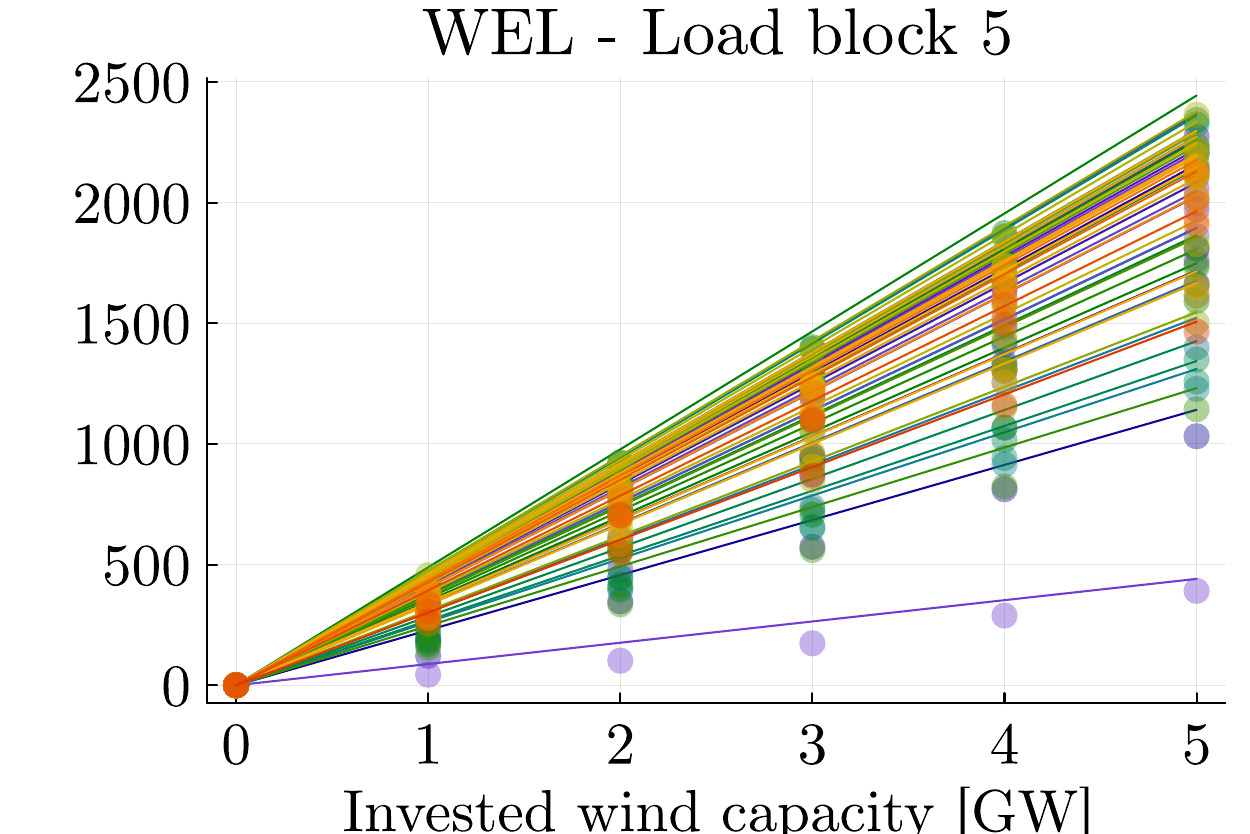}
\caption{Linear wind representation in Wellington for load block 5}
\end{subfigure}
\caption{Wellington wind generation (MW) versus $K$ (GW), where load block indices are determined using $\bar{K}$=2.5 GW. Different colours correspond to different weeks of the year.} 
\label{fig:windwel}
\end{figure}

The approximation can be illustrated by showing how it applies to some regions in New Zealand, as a preview to the next section where we discuss the New Zealand system in more detail. Figure~\ref{fig:windotg} shows the fitted linear functions representing wind generation for the Otago region of New Zealand for the first (peak), third and fifth load block.
A linear regression gives a good representation of the wind generation in each load block at least at the five data points we have chosen. The figure also shows that for most weeks, more of the wind generation is assigned to the load blocks with less demand (the average generation across the weeks shifts upwards as the block number increases). In fact, for all the regions the yearly mean wind generation is always increasing going from the first block to block two, three, four and five. 
Figure~\ref{fig:windwel} shows the data points and the corresponding linear regression for the Wellington region. This has a poorer fit, and a few of the  estimated slopes in load block 1 are negative. They end up negative because $W$ is discontinuous in $K$ (possibly) jumping downwards when $b(K)$ changes. The impact of this is strongest in the Wellington region because the value of $K_r/K$ is 
large 
compared to other regions in the model.

\section{A New Zealand case study}\label{sec:case}

In this section, we show the results of applying our model to the New Zealand electricity system. 
We do this by amending JADE, an existing open-source SDDP model of the New Zealand electricity system \citep{jadejl}. JADE is written in Julia \citep{bezanson_julia_2017} using the JuMP \citep{Lubin2023} and SDDP.jl \citep{dowson2021sddp} packages, and it is distributed by the New Zealand Electricity Authority. For this paper, we extended JADE to add the investment decisions and a wind generation model as described in the previous sections.

\begin{figure}[htbp]
\centering
\includegraphics[width=\textwidth]{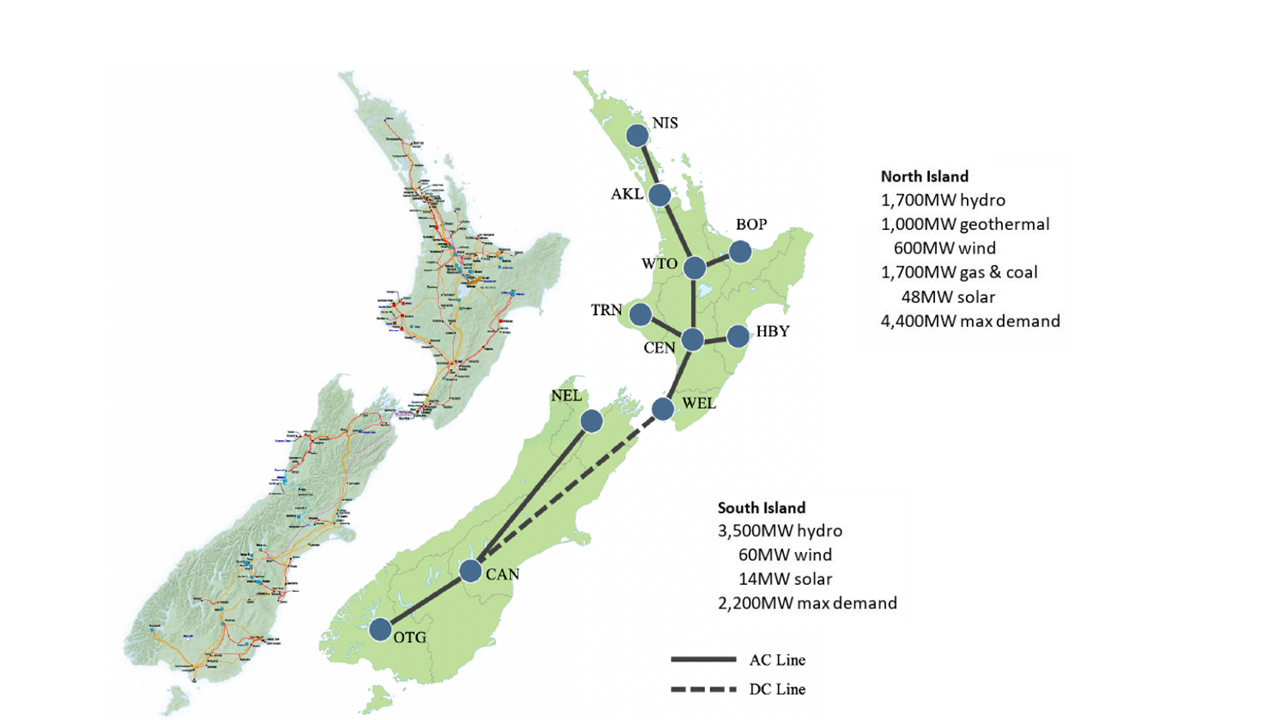}
\caption{The New Zealand electricity system, showing major transmission lines and the 11 node approximation used in the JADE model.} 
\label{fig:nzgrid}
\end{figure}

The New Zealand electricity system is spread over two islands as shown in Figure~\ref{fig:nzgrid}. The details of this system can be found in, for example, \cite{philpott2019new} and so we only give a brief description here. As in  \cite{philpott2023onslow}, JADE in the current paper approximates the full 250-node transmission network by an 11-node model as shown in the right-hand panel. Generation capacities and demand in each island for 2017 are also shown in Figure~\ref{fig:nzgrid}. The islands are joined by a HVDC cable of capacity 1050 MW (shown by the dashed line in Figure~\ref{fig:nzgrid}).

JADE has the ability to model arbitrarily many state variables but following \cite{philpott2023onslow} we limit these to storage volumes in Lake Taupo at node WTO, storage volumes in Lakes Tekapo, Ohau and Pukaki at node CAN, and storage volumes in Lake Hawea and Lakes Manapouri-Te Anau at node OTG. The proposed Lake Onslow reservoir is located in node OTG. The details of the power plants associated with these reservoirs are given in  \ref{appendix:system}, and the mathematical formulation of JADE is provided in \ref{appendix:jade}.

Our use of JADE enables a straightforward comparison to be made with the wind investment decisions optimized by \citet{philpott2023onslow} using enumeration in JADE. They considered three alternative scenarios for a fully renewable electricity system in 2035, in which all capacity decisions were determined up front apart from wind capacity that was optimized. In all three scenarios existing non-renewable electricity plants were closed. In the first case a large pumped hydro storage facility \textit{Lake Onslow} (1.5 GW and 5 TWh), is built to balance the system. The second alternative includes a zero-emission green peaker plant, and the third case assumes existing renewable capacity. We refer to these three cases as \textit{Onslow}, \textit{Peakers} and \textit{Wind only} from now on. The investments in new wind capacity for each case were then found by solving  JADE a number of times with different capacity levels and selecting the capacity choice yielding a generation weighted average price (GWAP) equal to the LCOE for the new wind power plants. The results of these experiments provide a useful benchmark for the investments generated by our method. 

Of the 11 regions shown in Figure~\ref{fig:nzgrid}, seven  have consented wind farms. We allocate new wind capacity proportional to the existing capacity in these regions, as shown in Table~\ref{tab:windfarms}. This means that the state variable $x_{inv}$ is implemented as the national wind capacity $K$, and the regional shares of this capacity are fixed to $\alpha_rK$. 

\begin{table}[!ht]
\centering
\begin{tabular}{c c c c}
Region & $\alpha_r$ & Representative wind power plant & Capacity [MW] \\
\hline
CAN & 0.0814 & White Hill & 58  \\
CEN & 0.0627 & Tararua Stage 3 & 93   \\
HBY & 0.0914 & Te Uku &  28 \\
OTG & 0.1926 & White Hill & 58 \\
TRN & 0.000 & Waipipi & 133 \\
WEL & 0.5670 & West Wind & 142.6 \\ 
WTO & 0.0049 & Te Uku & 28
\end{tabular}
\caption{Representative wind farms and $\alpha_r$ for each region with wind.}
\label{tab:windfarms}
\end{table}

Recall that the levelized cost of energy (LCOE) is defined as:
$$\text{LCOE} = \frac{\text{Lifetime cost of capacity} }{\text{Lifetime energy produced}}.$$
If we denote the overnight capital cost of 1 MW of wind by $I$, its lifetime by $\tau$ years, its capacity factor by $\eta$, and $H$ as the number of hours in a year (8760), then (ignoring maintenance costs):
$$\text{LCOE} = \frac{I}{H \cdot \eta \cdot ( 1 + \beta + \beta^2 + \ldots \beta^{\tau-1})}.$$
Thus:
\begin{equation*}
I = \text{LCOE} \cdot H \cdot \eta \cdot \frac{1-\beta^{\tau}}{1-\beta}.
\end{equation*}

Following \cite{philpott2023onslow}, we assume LCOE = 65 NZD/MWh. The national wind capacity factor $\eta$ is estimated from the wind profiles to be $0.355$, and we assume $\tau=20.$ When $\beta = 0.9$, this gives an overnight investment cost for our infinite-horizon model of $I~=~1.78\times 10^6$~NZD/MW.

\section{Results}
\label{sec:results}

In this section, we present the results of applying our model to the New Zealand case study briefly introduced in the previous section. The main results presented in this paper are the resulting invested wind capacities for different model runs with our altered model. 

\subsection{Comparison to Philpott \& Downward (2023)}

\begin{figure}[!htbp]
\centering
\includegraphics[width=0.7\textwidth]{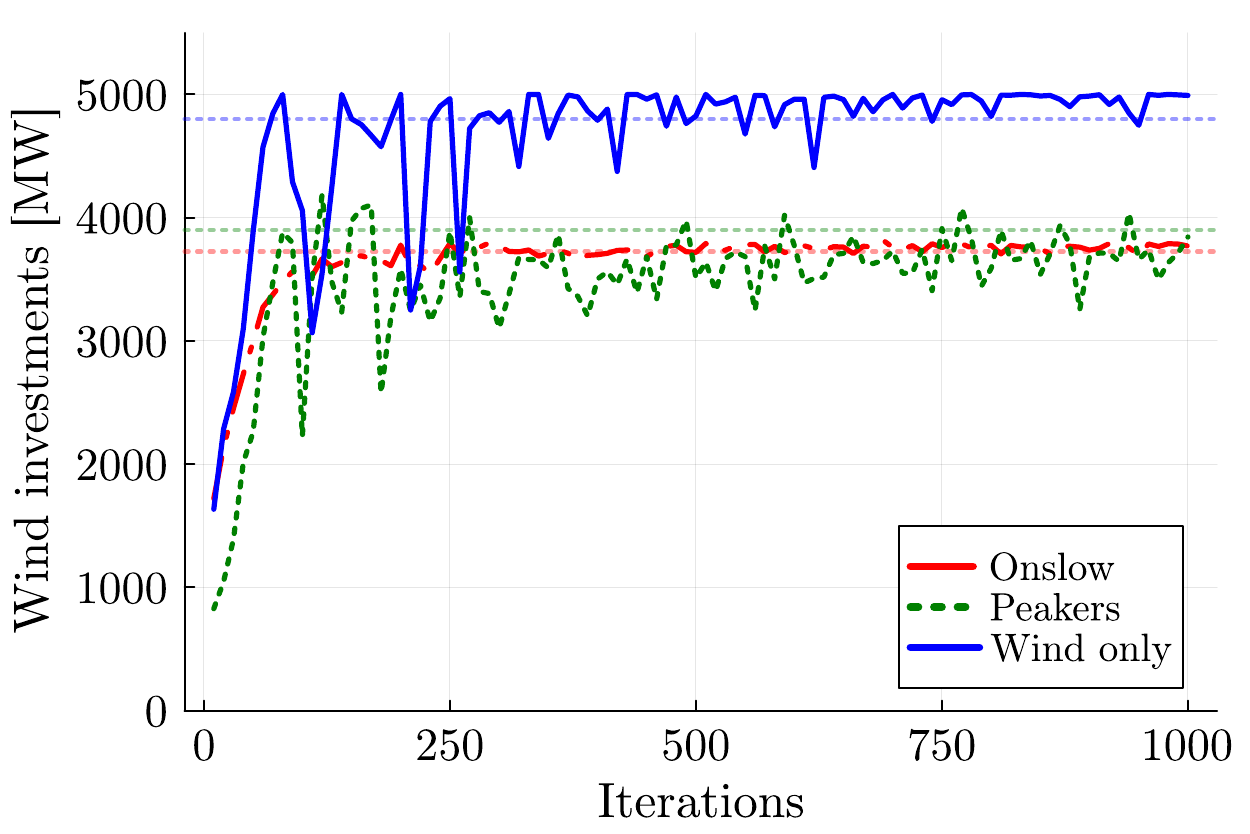}
\caption{Investment decisions for new wind capacity in the three different cases against the number of SDDP training iterations. Horizontal, dotted lines are those found by \cite{philpott2023onslow}.} 
\label{fig:comparinginvestments}
\end{figure}

Figure~\ref{fig:comparinginvestments} compares the investment decisions from our model with the levels published in \cite{philpott2023onslow}, which were found by re-solving JADE with different investment levels and choosing the best (shown by the horizontal, dotted lines). In each of the three cases, our model is trained for ten SDDP iterations, and then the investment decision for the national wind capacity is evaluated. This is repeated for a total of 1000 SDDP iterations, adding new cuts to the investment problem for ten consecutive forward and backward passes at the time. The plot shows that the endogenous investment decision approaches the benchmark investment levels as the model is trained with more iterations.

In general, the results in Figure~\ref{fig:comparinginvestments} indicates that our method for integrating the investment decision works, with resulting wind investments close to the ones found in the previous study. We would not expect the optimal investments in our JADE model to match exactly the optimal investments in the JADE model used in \cite{philpott2023onslow} for several reasons. First, the linear wind representation described in Section~\ref{sec:wind} is a simplification of the affect that wind investment will have on residual load, compared with \cite{philpott2023onslow} where the load duration curve is re-computed for each level of investment.  Second, our model represents tranches of load shedding as a proportion of original load, whereas \cite{philpott2023onslow} treats these as proportions of residual load (after wind has been subtracted). Finally, \cite{philpott2023onslow} uses round numbers when enumerating different levels of wind investment, whereas the investment decision in our model is a continuous variable.

\begin{figure}[htbp]
\centering
\begin{subfigure}[b]{0.31\textwidth}
\caption{}
\includegraphics[width=\textwidth]{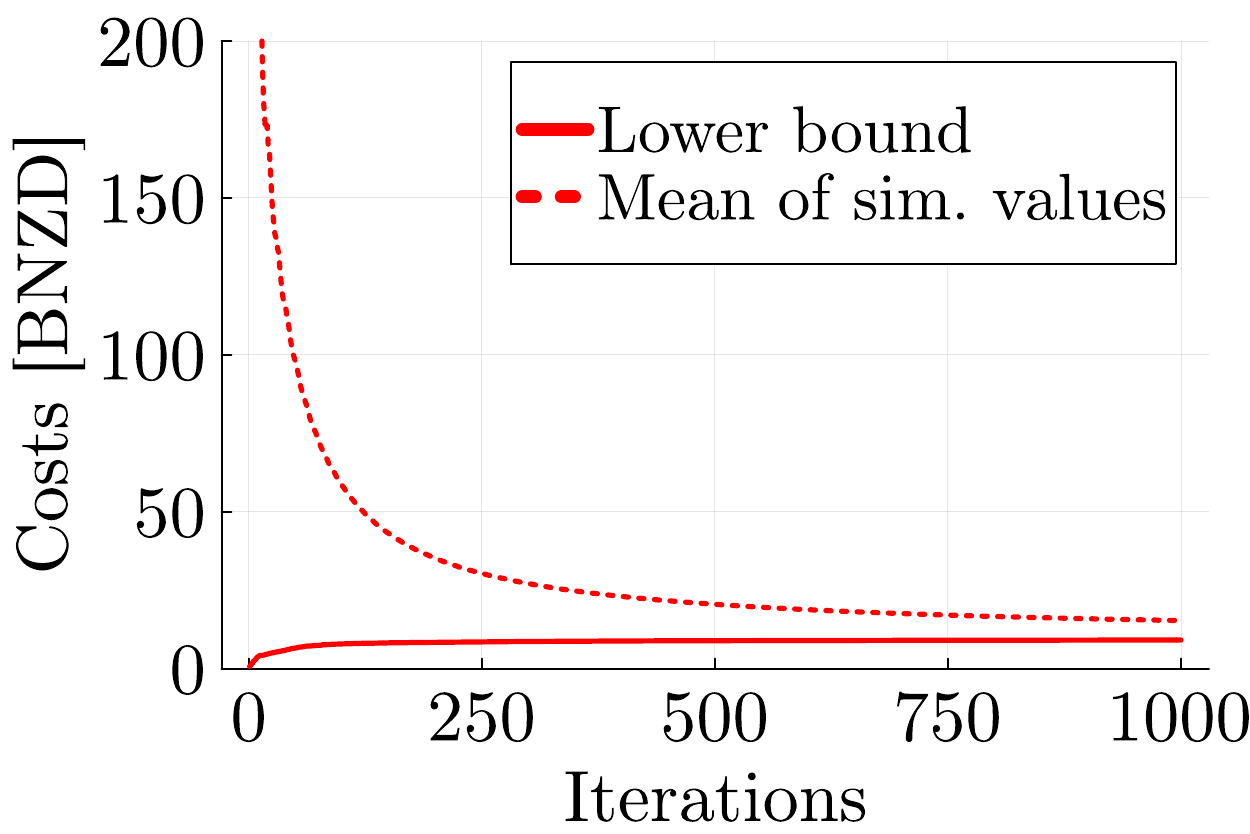}
\label{fig:convergenceonslow}
\end{subfigure} 
\hfill
\begin{subfigure}[b]{0.31\textwidth}
\caption{}
\includegraphics[width=\textwidth]{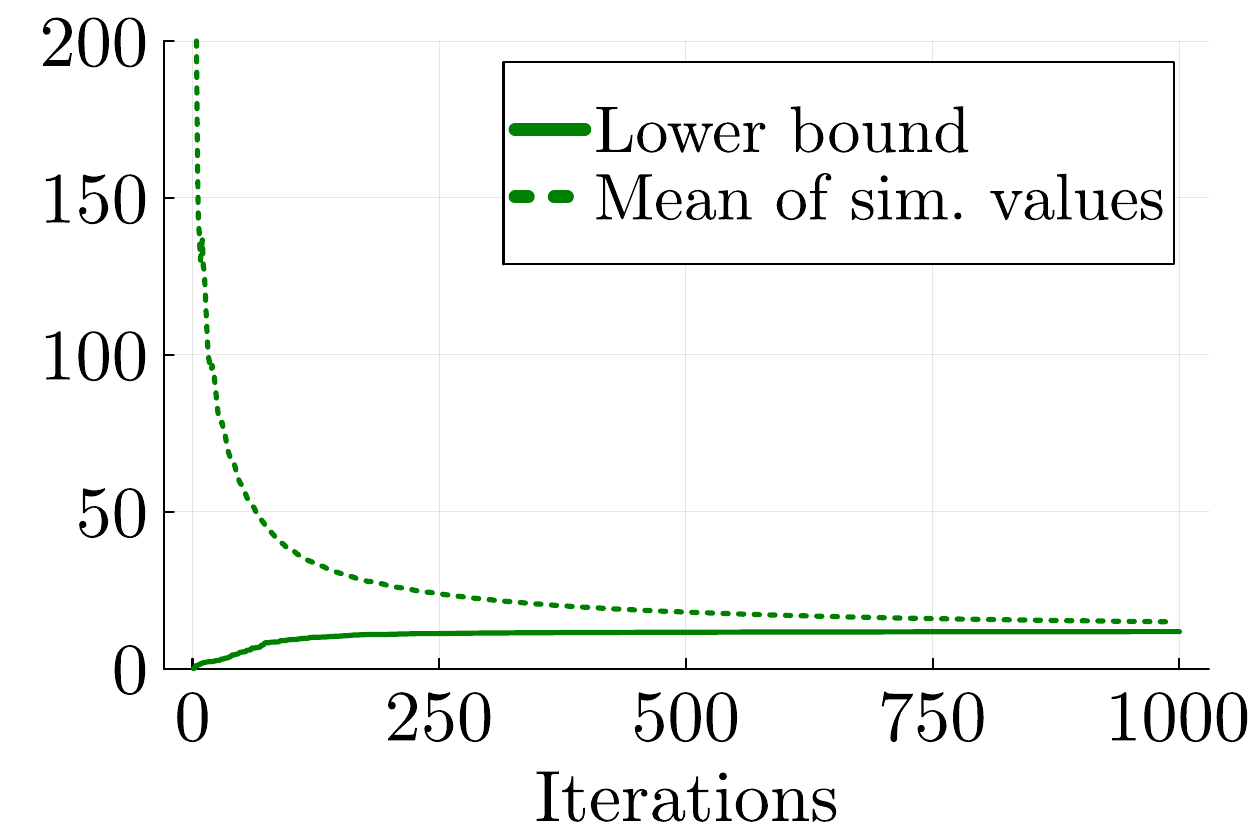}
\label{fig:convergencepeakers}
\end{subfigure} 
\hfill 
\begin{subfigure}[b]{0.31\textwidth} 
\caption{}
\includegraphics[width=\textwidth]{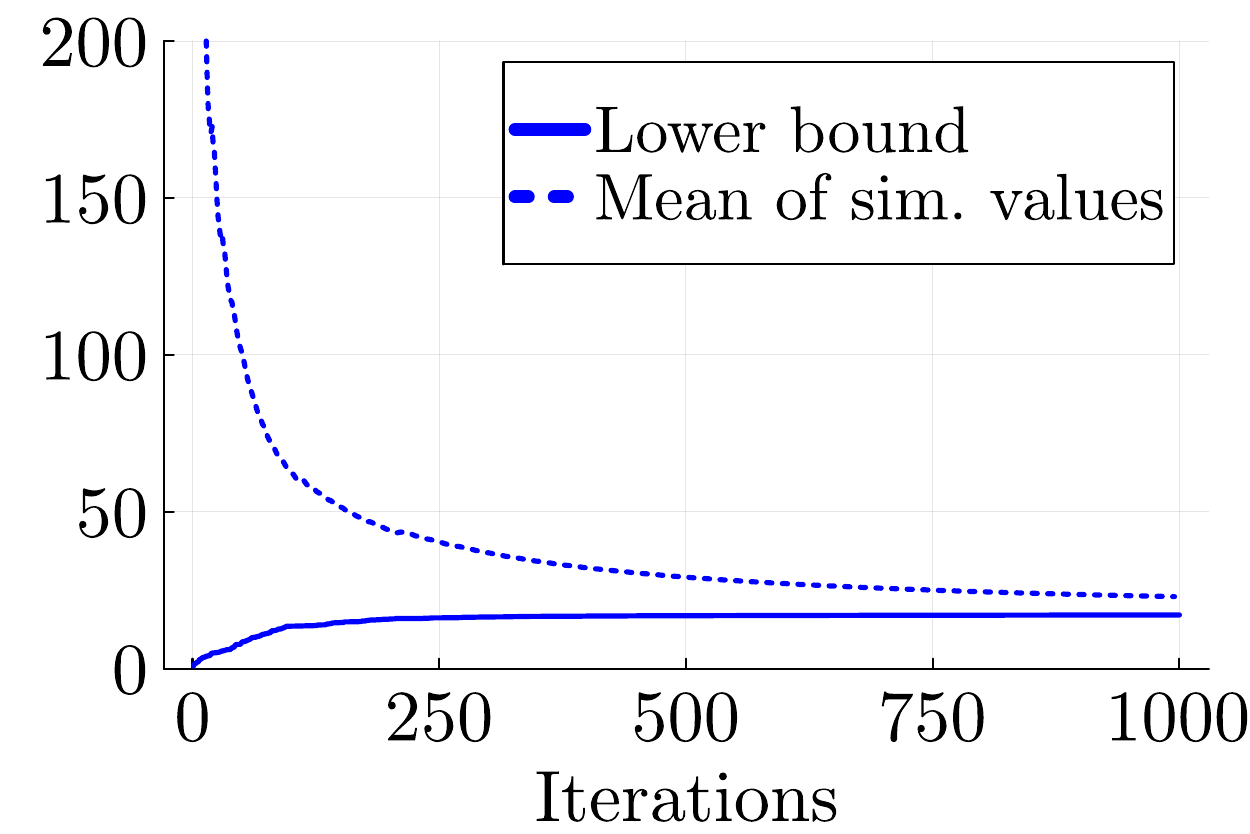} 
\label{fig:convergencewindonly} 
\end{subfigure}
\caption{Plots of the lower bound and the rolling mean of simulation values from the forward pass for each of the three cases when the model is trained with 1000 iterations: (a) \textit{Onslow}; (b) \textit{Peakers}; and (c) \textit{Wind only}} 
\label{fig:convergence}
\end{figure}

Despite some differences between the two versions of the model, the resulting investments are similar. Evaluating the policy every tenth iteration however, introduces some randomness to Figure~\ref{fig:comparinginvestments}. It is easy to see that the investments improve as more cuts are added to the policy, but the plot does not prove convergence in any way. Figure~\ref{fig:convergence} shows the lower bound and the rolling mean of the simulation values from the forward pass for all of the three cases. We terminated the training loop of SDDP.jl after a fixed iteration limit of 1000 iterations.

\subsection{Incorporating HVDC and peaker capacity}

One of the strengths of integrating the investment decisions in the model compared to the method from \cite{philpott2023onslow} is that one can easily optimize more investment types, at the cost of increased computational effort. To illustrate this, we conduct an experiment where the wind capacity is split into North Island and South Island investments. In addition, we let the model optimize the capacity of the green peaker in the \textit{Peakers} case, and the HVDC transmission capacity between the North Island and the South Island of New Zealand in all of the cases. Thus, the $x_{inv}$ state variable is now a vector with four elements, increasing the number of subproblems solved in training the policy with 1000 iterations from 11.6 million to 15.4 million. Solving such a case by enumeration is much more difficult.

\begin{table}[ht]
    \centering
    \begin{tabular}{r | c c c}
    Investments and costs & \textit{Onslow} & \textit{Peakers} & \textit{Wind only} \\
    \hline
    Wind - North Island [MW] & 1315 & 1458 & 1997  \\
    Wind - South Island [MW] & 0 & 0 & 0  \\
    HVDC capacity [MW] & 643 & 0 & 379 \\ 
    Green peaker [MW] & - & 837 & - \\
    \hline
    CAPEX [MNZD] & 4633 & 3671 & 5047 \\
    OPEX [MNZD] & 2128 & 2219 & 3775 \\
    Total [MNZD] & 6761 & 5890 & 8822
    
    \end{tabular}
\caption{Resulting investments and costs from the extended investment model. \textit{CAPEX} is the cost of the deterministic investment node. \textit{OPEX} is the SDDP lower bound less the CAPEX cost.}
\label{tab:inv4}
\end{table}

Table~\ref{tab:inv4} summarizes the results from the extended investment model after training with 1000 iterations. The total wind investment is lower in all cases than the values in Figure~\ref{fig:comparinginvestments}. Multiple factors are contributing to this outcome. In our previous model we assumed that the capacity on the HVDC connection between the South Island and the North Island was sufficient to meet all transfers. This enabled a comparison with  \cite{philpott2023onslow}  as shown in Figure~\ref{fig:comparinginvestments}.  The new model optimizes the HVDC capacity, with the existing capacity as the starting point. Thus, optimal wind investments are lower when wind investments needs to come with following grid investments to access storage in the South Island. In this case, shedding more load is cost efficient compared to extensive grid expansion. Also, as load centers are located on the North Island (and the capacity factor of the representative wind power plants on the North Island is higher than on the South Island), all wind investments are made in the North Island only.

Table~\ref{tab:inv4} shows that increased HVDC capacity is most beneficial in the \textit{Onslow} case, because the Lake Onslow pumped hydro storage is located in the South Island and requires grid capacity to provide its flexibility to the North Island.

In the \textit{Peakers} case, the optimal value of peaking capacity determined by the new model is significantly higher than that assumed in \cite{philpott2023onslow}. The \textit{Wind only} case has the highest CAPEX, because of the larger wind investments, and the highest OPEX, because of load shedding, of the three cases. The investment costs from building Lake Onslow (estimated in \cite{NZBattery} to  be 15700 MNZD) are not included in the CAPEX cost for the \textit{Onslow} case.  Thus the total cost for the \textit{Onslow} case is 22461 MNZD, which is almost four times the cost for the cheapest alternative, namely the \textit{Peakers} case.

\begin{figure}[!ht]
\centering
\includegraphics[width=0.7\textwidth]{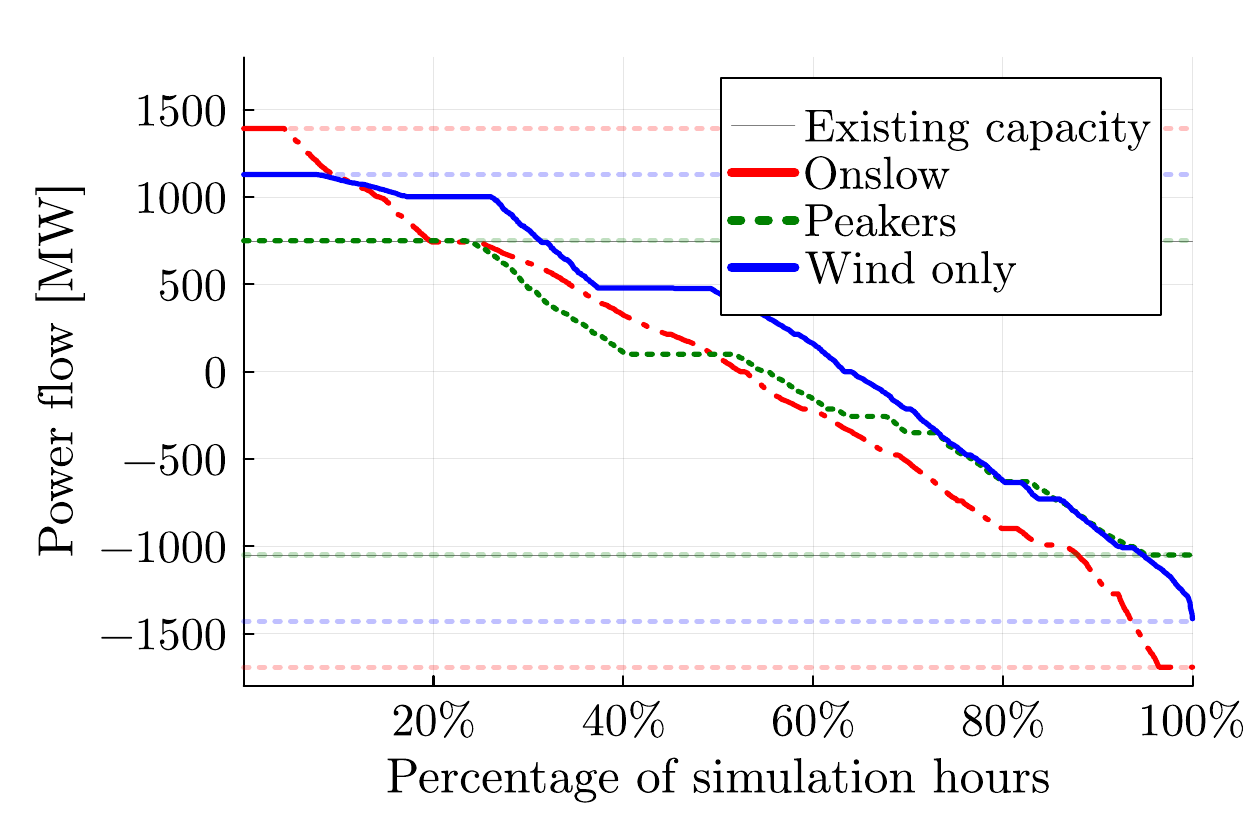}
\caption{Flow-duration curve for flow on the HVDC line from WEL (North Island) to CAN (South Island) in the three different cases when simulated over the 31 historical inflow years 1990--2020. The horizontal, dotted lines shows the maximum capacity of the lines in each case, which is the sum of the existing capacity (the gray line) and the invested grid capacity.}
\label{fig:comparinghvdcflow}
\end{figure}

Figure~\ref{fig:comparinghvdcflow} shows the utilization of the HVDC cable in each investment case over the 31 historical inflow years 1990--2020. The resulting HVDC capacities are shown as the horizontal, dotted lines, while the three plotted lines display the duration curve for the load flow on the HVDC connection with positive load flow defined as flow from the Wellington region (North Island) to the Canterbury region (South Island). The green, horizontal line is equal to the existing capacity (gray line) as now grid investments are made in the \textit{Peakers} case. The $x$-axis gives the percentage of simulated hours where load flow is equal to or greater than $y$ MW. In all cases the  utilization of the transmission line goes both ways, with the largest flows in the \textit{Onslow} case, where North to South flow is needed for pumping and reverse flow transports Onslow generation north. In the two other cases, HVDC investments are lower, and the full capacity from the North Island to the South Island is utilized more that the full reverse flow.

\begin{figure}[!ht]
\centering
\includegraphics[width=0.8\textwidth]{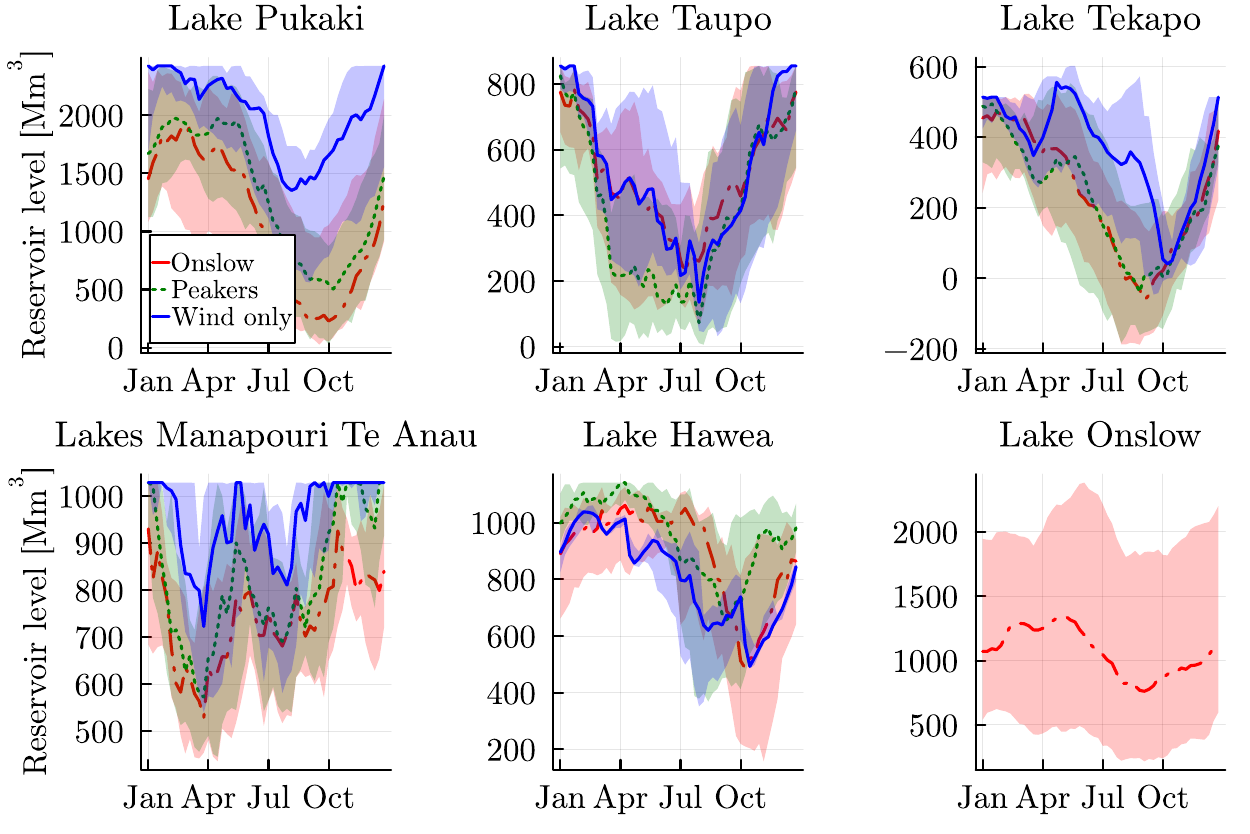}
\caption{Hydro reservoir levels relative to consented minimum operating levels for the three different cases when simulated over the 31 historical inflow years 1990--2020. The lines are the median reservoir level and each coloured band depicts the 25th and 75th percentile levels. Some reservoirs can have negative storage when the reservoir level falls below the consented operating minimum (incurring a cost penalty).} 
\label{fig:comparingstorage}
\end{figure}

Figure~\ref{fig:comparinghvdcflow} is complemented by Figure~\ref{fig:comparingstorage}, which visualizes the utilization of the six (five without Lake Onslow) main reservoirs in the system. The figure shows that four of the five reservoirs are operated with lower storage levels when other sources of flexibility are in place, like Lake Onslow or a green peaker plant, hence reducing the risk of spillage.

\begin{figure}[!ht]
\centering
\includegraphics[width=0.8\textwidth]{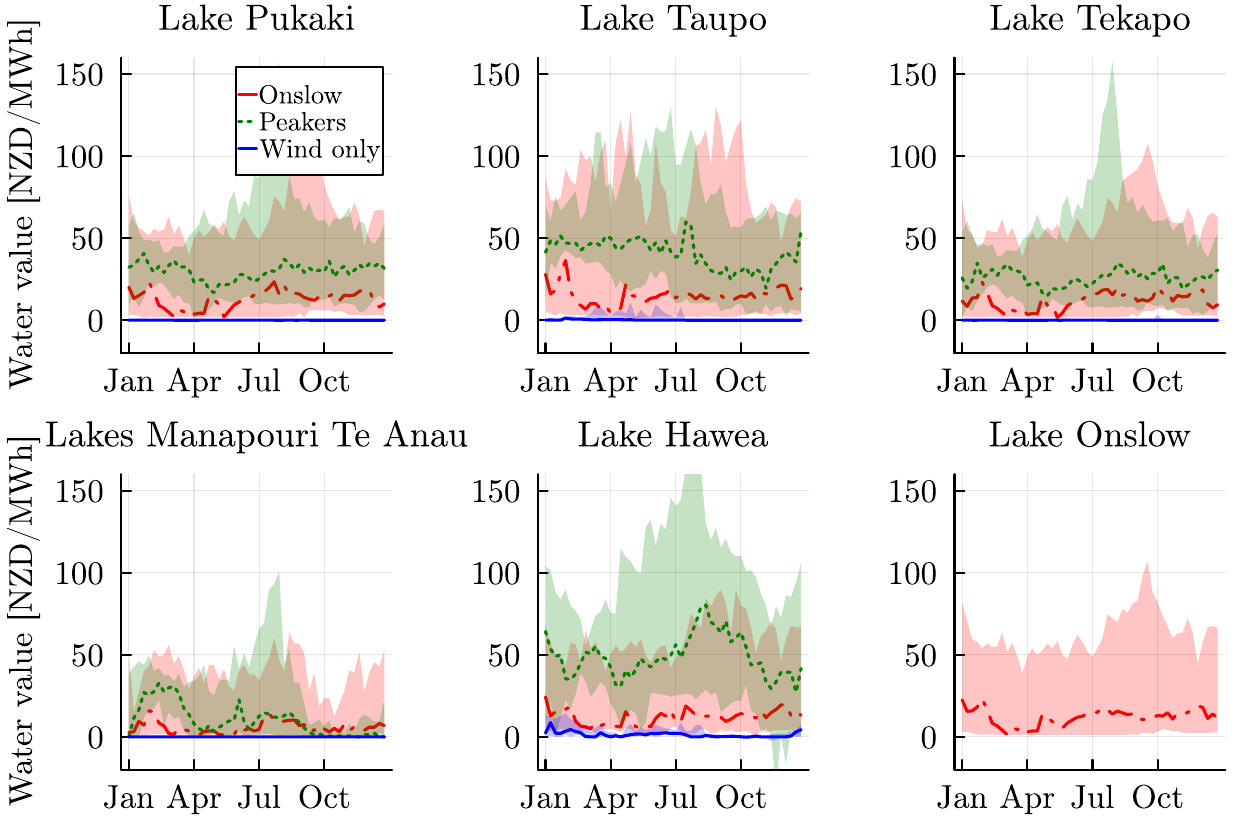}
\caption{Marginal water values for different reservoirs in the three different cases when simulated over the 31 historical inflow years 1990--2020. The lines are the median marginal water value and each coloured band depicts the 25th and 75th percentiles.} 
\label{fig:comparingwatervalues}
\end{figure}

Figure~\ref{fig:comparingwatervalues} plots trajectories of marginal water values in the hydro reservoirs. These values can be taken as an indication of spot energy prices when these are determined by energy  constraints rather than capacity constraints.  Marginal water values are very low in the {\it Wind only} case, and higher in the other two cases. In the {\it Wind only} case, prices are set by periods when the wind is not blowing and load must be shed because conventional capacity and/or grid capacity is insufficient, even though energy in reservoir storage is plentiful giving low marginal water values. In the other two cases, we have increased dispatchable capacity, and prices are set by anticipated energy shortages that are reflected by higher marginal water values.
More detailed plots of the storage trajectories and marginal water values for the 31 historical inflow years from 1990--2020 can be found in \ref{appendix:storage} and \ref{appendix:watervalue}. 

\begin{figure}[!ht]
\centering
\includegraphics[width=0.7\textwidth]{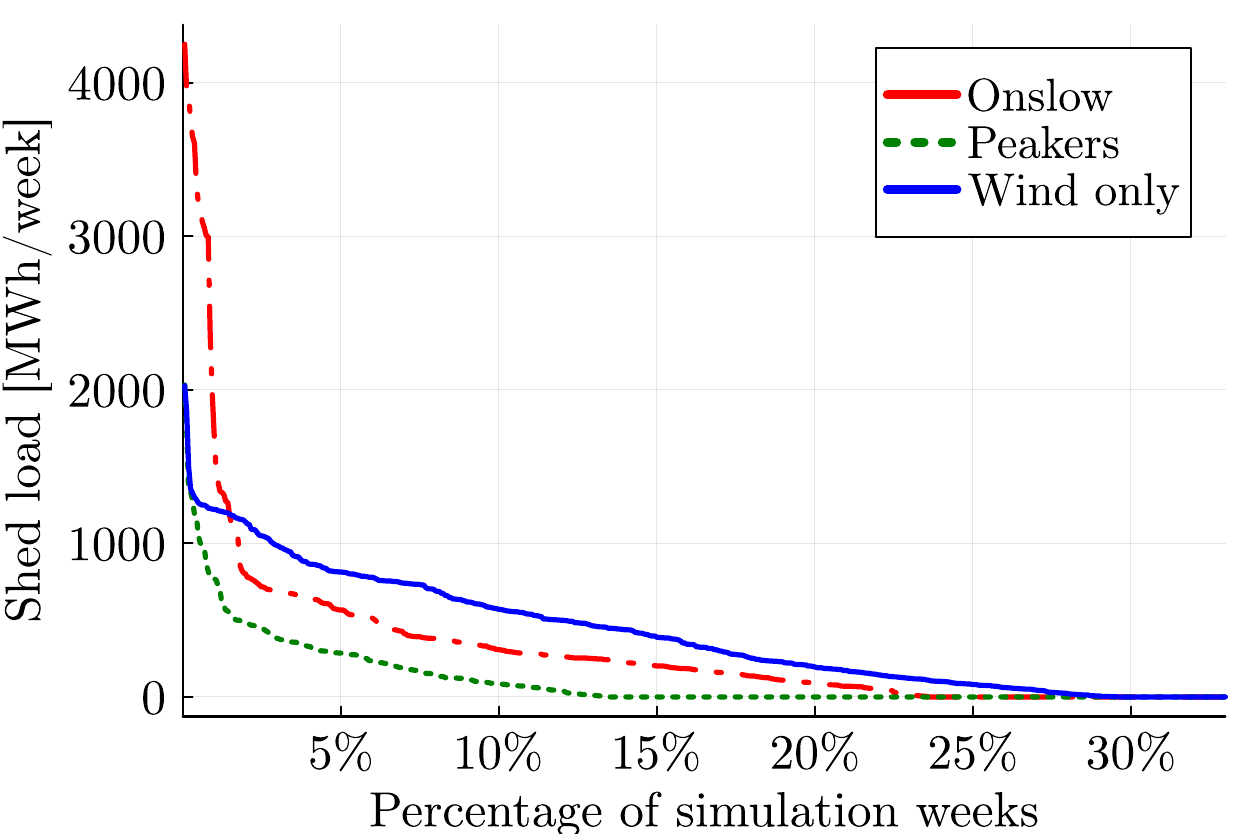}
\caption{Duration curves of the weekly load shedding in the three cases when simulated over the 31 historical inflow years 1990--2020.} 
\label{fig:shedding}
\end{figure}

When demand cannot be met, the system must shed load. Figure~\ref{fig:shedding} shows the probability of shedding more than $y$ MWh of load in a random week. In all three cases, more than 10 percent of the 1612 simulated weeks have some load shedding.

The \textit{Wind only} case sheds load in a large proportion of weeks, in periods when the amount of wind generation is unable to meet demand even if all other generators are at full capacity. This means prices are set by shortage costs while marginal water values are generally lower (since generation or grid capacity constraints mean extra water cannot be used to offset this shortage). The prices in these periods provide revenue for wind investment; in other periods prices are zero.

The \textit{Peakers} case provides dispatchable generation capacity to reduce the number of shortage periods.  The green trajectory in Figure~\ref{fig:shedding} shows that this results in less load shedding from capacity constraints. Some of the investment that was made in the \textit{Wind only} case is diverted to peakers. The availability of more dispatchable plant makes reservoir operation less conservative making energy shortages more likely, and marginal water values increase.

The \textit{Onslow} case also provides dispatchable generation capacity to reduce the number of shortage periods so the red trajectory in Figure~\ref{fig:shedding} lies below the blue, at least for nearly every week in the simulation. Observe that there are a small number of weeks in the simulations when the load shedding is much higher. As observed in \cite{philpott2023onslow}, these cases occur when Onslow generation is required but the lake has been emptied by a dry period in the previous year. Marginal water values and load shedding become very high in these weeks.

\section{Conclusions}
\label{sec:conclusions}

In this paper, we have shown how capacity expansion investments can be included in an infinite-horizon operational SDDP model by augmenting the policy graph to include an investment node. We derived a linear representation of wind generation in different load blocks to link the investment decisions with the subproblem constraints. Additionally, by applying our method to the New Zealand energy system, we have shown that our proposed method is flexible, and that it can be used to analyse a variety of options for expansion within the same system. Because our model is based on the policy graph and SDDP.jl, it did not require us to code a customized decomposition algorithm, which makes our model easier to adopt and implement than alternatives in the literature.

The results from our case study show that, with the assumptions we made, a peaker plant located in the North Island is a more cost efficient supplement to wind investments than the \textit{Lake Onslow} pumped hydro storage.\footnote{As we were finishing writing this article, the New Zealand Government announced that they will scrap the plans to build Lake Onslow \citep{scraponslow}.}

\section*{Acknonwledgements}

J. Hole was supported by the Norwegian Water Resources and Energy Directorate (NVE), the Department of Electric Energy at the Norwegian University of Science and Technology (NTNU) and the NZ MBIE Catalyst grant UOCX2117 - New Zealand-German Platform for Green Hydrogen Integration (HINT). The work is associated with the two Centres for Environment-friendly Energy Research FME CINELDI (RCN grant 257626) and FME NTRANS (RCN grant 296205).

\bibliography{references}

\appendix

\section{JADE formulation}\label{appendix:jade}

JADE seeks a policy of electricity generation that meets demand and minimizes the expected cost of thermal generation fuel consumed plus any costs of load reduction. All data are deterministic except for weekly inflows that are assumed to be stagewise independent. The resulting stochastic dynamic programming model is defined as follows. Let $x_{j}\left( t\right) $ denote the storage in reservoir $j$ at the \emph{end} of week $t$, and let the Bellman function $C_{t}(\bar x,\omega (t))$ be the minimum discounted expected fuel cost to meet electricity demand in weeks $t,t+1,\ldots$, when reservoir storage $x_{j}(t-1)$ at the start of week $t$ is equal to $\bar{x}_{j}$ and the inflow to reservoir $j$ in week $t$ is known to be $\omega_{j}(t)$. 

In JADE, a weekly discount factor $\rho=\beta^{1/52}$ is used when going from any stage to the next, where $\beta < 1$ is the annual discount factor. This implies that in each stage there is a probability of $1-\rho$  of transitioning from stage t to a zero-node $0$ where $C_{0}(\bar x,\omega (0)) = 0$. 

The Bellman function $C_{t}(\bar x,\omega (t))$ for week $t$ is the optimal solution value of the mathematical program:

\begin{equation*}
\begin{array}{rll}
C_{t}(\bar{x},\omega(t)) = \\\min & \sum_{i\in \mathcal{N}} \sum_{b}T(b,t) \left( \sum_{m\in \mathcal{F}(i)}\phi _{m} f_{m}(b,t) 
+   \sum_{l\in \mathcal{L}(i)} \psi _{lb} z_{i}(l,b,t) \right) \\
& + \;\rho \cdot \mathbb{E}[C_{t+1}(x(t),\omega (t+1))]\\
\text{s.t.} & g_{i}(y(b,t))+\sum_{m\in \mathcal{F}(i)}f_{m}(b,t)\; + \\
& \qquad \sum_{m\in \mathcal{H}(i)}\gamma _{m}h_{m}(b,t) + \sum_{l\in \mathcal{L}(i)} z_{i}(l,b,t) =D_{i}(b,t), & i\in 
\mathcal{N}\\
& x(t)=\bar{x}-S\sum_{b}T(b,t)\left(A\,h(b,t)+A\,s(b,t)-\omega (t)\right)\\
& 0\leq f_{m}(t)\leq a_{m}, & m\in \mathcal{F}(i),\ i\in 
\mathcal{N}\\
& 0\leq h_{m}(t)\leq b_{m},\quad 0\leq s_{m}(t)\leq c_{m}, & m\in \mathcal{H}(i)\\
& 0\leq x_{j}(t)\leq r_{j}, & j\in \mathcal{J},\ i\in 
\mathcal{N},\ y\in Y.
\end{array}
\end{equation*}

This description uses the following indices:%
\begin{equation*}
\begin{array}{ll}
\text{Index} & \text{Refers to} \\ 
t & \text{index of week} \\ 
i & \text{node in transmission network} \\ 
b & \text{index of load block} \\ 
m & \text{index of plant} \\ 
j & \text{index of reservoir} \\ 
\mathcal{N} & \text{set of nodes in transmission network} \\ 
\mathcal{F}(i) & \text{set of green peaker plants at node }i \\ 
\mathcal{H}(i) & \text{set of hydro plants at node }i \\
\mathcal{L}(i) & \text{set of load types at node }i \\ 
\mathcal{J} & \text{set of reservoirs.}%
\end{array}%
\end{equation*}%

The parameters are:%
\begin{equation*}
\begin{array}{lll}
\text{Symbol} & \text{Meaning} & \text{Units} \\ 
\phi _{m} & \text{short-run marginal cost of peaker plant }m & \text{\$/MWh}
\\ 
\psi _{lb} & \text{cost of shedding load type }l \text{ in load block }b & \text{\$/MWh}
\\
\gamma _{m} & \text{conversion factor for water flow into energy} & \text{%
MWs/m}^{3} \\ 
D_{i}(b,t) & \text{electricity demand in node }i\text{ in block }b\text{,
week }t & \text{MW} \\ 
T(b,t) & \text{number of hours in load block }b\text{ in week }t & \text{h}
\\ 
S & \text{number of seconds per hour (3600)} &  \\ 
a_{m} & \text{thermal plant capacity} & \text{MW} \\ 
b_{m} & \text{hydro plant capacity} & \text{m}^{3}\text{/s} \\ 
c_{m} & \text{spillway capacity} & \text{m}^{3}\text{/s} \\ 
r_{j} & \text{reservoir capacity} & \text{m}^{3} \\ 
Y & \text{feasible set of transmission flows} &  \\ 
A & \text{incidence matrix of river chain} & 
\end{array}%
\end{equation*}

The variables are:
\begin{equation*}
\begin{array}{lll}
\text{Symbol} & \text{Meaning} & \text{Units} \\ 
\rho & \text{weekly discount factor} & \\
f_{m}(b,t) & \text{generation of green peaker plant }m\text{ in load block }b%
\text{ in week }t & \text{MW} \\ 
z_{i}(l,b,t) & \text{shed load of type }l\text{ in load block }b \text{ in node }i
\text{ in week }t & \text{MW} \\ 
x_{j}(t) & \text{storage in reservoir }j\text{ at end of week }t & \text{m}%
^{3} \\ 
\bar{x}_{j} & \text{known storage in reservoir }j\text{ at start of week }t
& \text{m}^{3} \\ 
h(b,t) & \text{vector of hydro releases in block }b\text{, week }t & \text{m}%
^{3}\text{/s} \\ 
s(b,t) & \text{vector of hydro spills in block }b\text{, week }t & \text{m}%
^{3}\text{/s} \\ 
\omega (t) & \text{inflow (assumed constant over the week)} & \text{m}^{3}%
\text{/s} \\ 
y(b,t) & \text{flow in transmission lines in load block }b\text{ in week }t
& \text{MW} \\ 
g_{i}(y) & \text{sum of flow into node }i\text{ when transmission flows are }%
y & \text{MW} \\ 
\end{array}%
\end{equation*}

The water-balance constraints in the storage reservoirs at the end of week $t$ are represented by:
\begin{equation*}
x(t)=\bar{x}-S\sum_{b}T(b,t)(A\,h(b,t)+A\,s(b,t)-\omega (t)),
\end{equation*}%
where $x_{j}(t)$ is the storage in reservoir $j$ at the end of week $t$, $%
s_{j}(b,t)$ denotes the rate of spill (in m$^{3}$/second) in load block $b$
in week $t$, and $\omega _{j}(t)$ is the uncontrolled rate of inflow into
reservoir $j$ in week $t$. We multiply all of these by $S$ to convert to m$%
^{3}$/hour, and then by $T(b,t)$ to give m$^{3}$ in each load block. All
these are subject to capacity constraints. (In some cases we also have
minimum flow constraints that are imposed by environmental resource
consents.) The parameter $\gamma _{m}$, which varies by generating station $%
m $, converts flows of water $h_{m}(t)$ into electric power. The same variables and constraints can be used to model pumping of water into a reservoir, except the value of the parameter $\gamma_m$ is negative to reflect that energy is consumed as water is pumped into a higher reservoir.

\section{System data}\label{appendix:system}
We summarize here the features of the New Zealand system represented in JADE. A full data set can be downloaded from  \cite{EAjade}. For our study, all existing fossil-fuel plant have been removed from the model.

\begin{table}[h]
\centering
\begin{adjustbox}{max width=215pt}
\begin{tabular}{c c r r}
Generator & Region & Capacity [MW]  & Specific power (MW/cumec) \\
\hline
Arapuni	&	WTO	&	192	&	0.462	\\
Aratiatia	&	WTO	&	78	&	0.284	\\
Atiamuri	&	WTO	&	84	&	0.196	\\
Karapiro	&	WTO	&	96	&	0.264	\\
Maraetai	&	WTO	&	352	&	0.526	\\
Matahina	&	BOP	&	80	&	0.595	\\
Ohakuri	&	WTO	&	112	&	0.284	\\
Rangipo	&	CEN	&	120	&	1.960	\\
Tokaanu	&	CEN	&	240	&	1.750	\\
Waikaremoana	&	HBY	&	140	&	3.535	\\
Waipapa	&	WTO	&	54	&	0.139	\\
Whakamaru	&	WTO	&	124	&	0.316	
\end{tabular}
\end{adjustbox}
\caption{North Island hydro generation stations optimized in JADE model}
\label{tab:NIgenerators}
\end{table}

\begin{table}[h]
\centering
\begin{adjustbox}{max width=215pt}
\begin{tabular}{c c r r}
Generator & Region & Capacity [MW]  & Specific power (MW/cumec) \\
\hline
Aviemore	&	CAN	&	220	&	0.310	\\
Benmore	&	CAN	&	540	&	0.818	\\
Clyde	&	OTG	&	464	&	0.518	\\
Cobb	&	NEL	&	32	&	4.405	\\
Coleridge	&	CAN	&	39	&	1.009	\\
Manapouri	&	OTG	&	842	&	1.531	\\
Ohau\_A	&	CAN	&	264	&	0.501	\\
Ohau\_B	&	CAN	&	212	&	0.417	\\
Ohau\_C	&	CAN	&	212	&	0.417	\\
Roxburgh	&	OTG	&	320	&	0.402	\\
Tekapo\_A	&	CAN	&	27	&	0.232	\\
Tekapo\_B	&	CAN	&	154	&	1.285	\\
Waitaki	&	CAN	&	105	&	0.162	\\
Onslow\_Pump	&	OTG	&	1500	&	-7.027	\\
Onslow\_Gen	&	OTG	&	1500	&	5.417	
\end{tabular}
\end{adjustbox}
\caption{South Island hydro generation stations optimized in JADE model}
\label{tab:SIgenerators}
\end{table}

\begin{figure}[htbp]
\centering
\includegraphics[width=\textwidth]{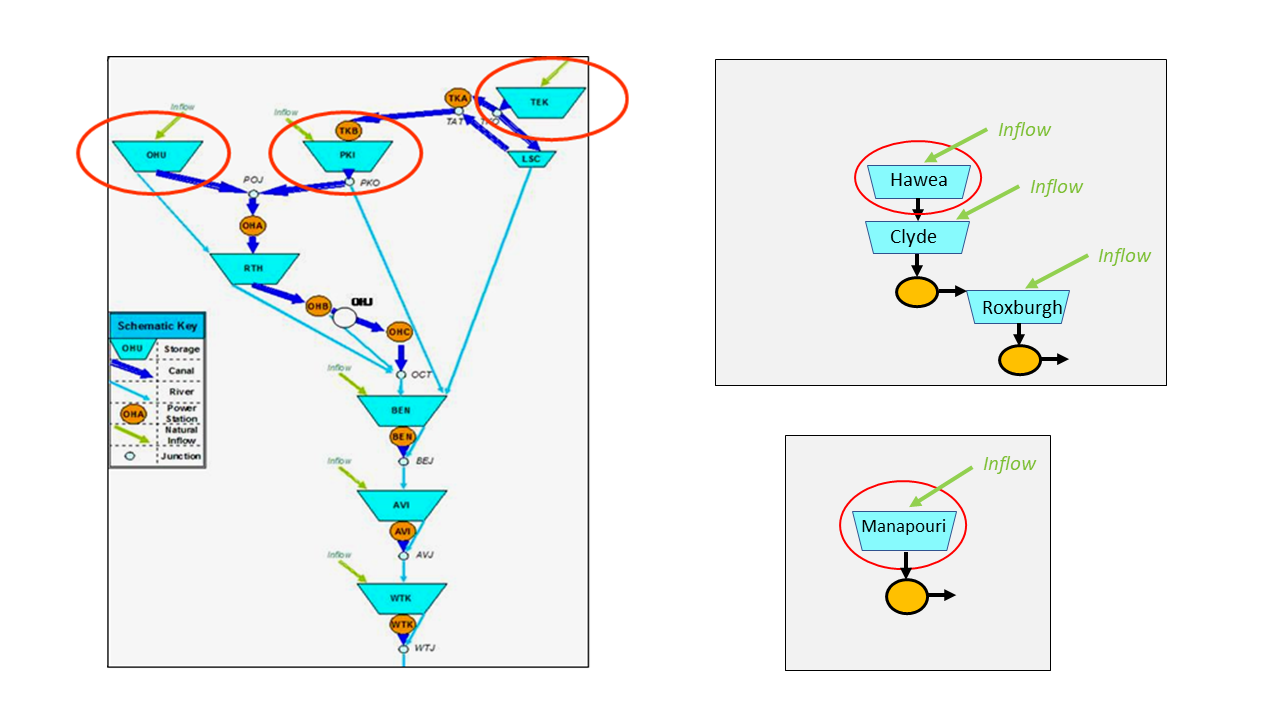}
\caption{South Island river chains optimized in JADE. Storage in circled reservoirs are state variables. South Island hydro plant not shown here are assumed to be run-of-river.
}
\label{fig:SIChains}
\end{figure}

\begin{figure}[htbp]
\centering
\includegraphics[width=\textwidth]{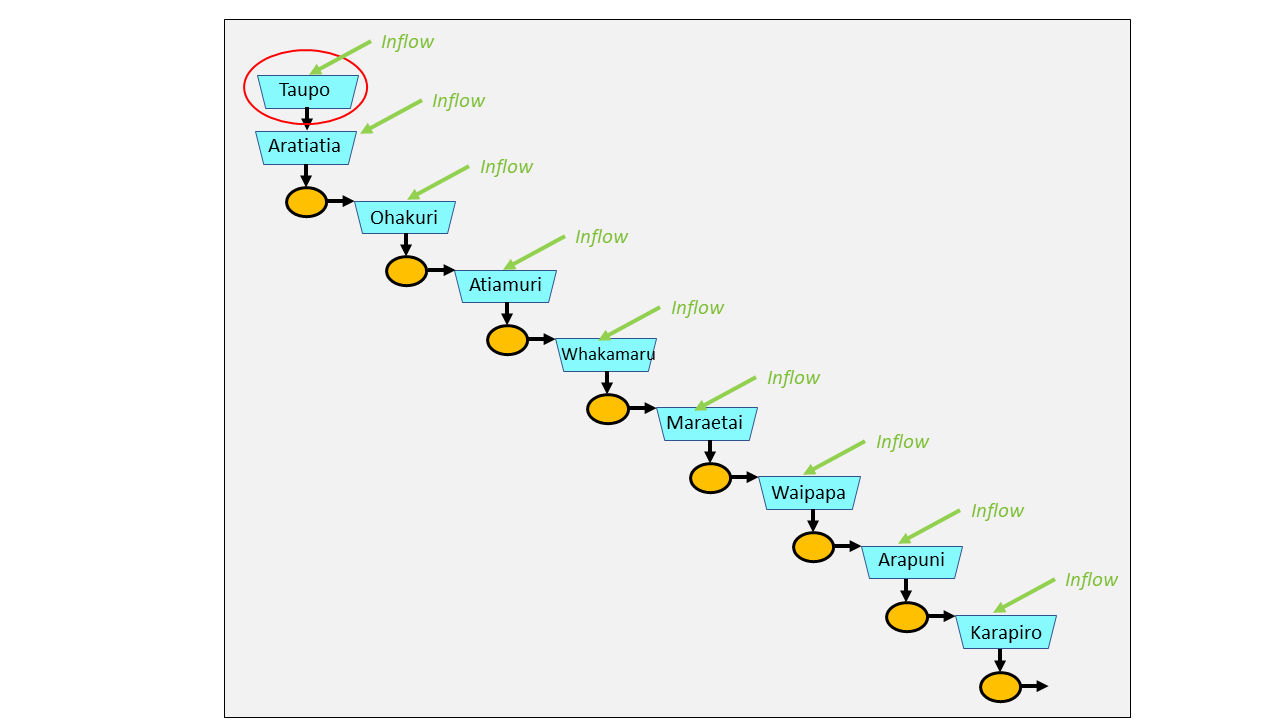}
\caption{North Island river chain (Waikato river) optimized in JADE. Storage in circled reservoir (Taupo) is a state variable. North Island hydro plant not shown here are assumed to be run-of-river.
}
\label{fig:NIChain}
\end{figure}

\begin{table}[h]
\centering
\begin{adjustbox}{max width=215pt}
\begin{tabular}{c r}
Region & Generation (MW) \\
\hline
AKL	&	39.54	\\
WTO	&	864.97	\\
BOP	&	240.23	\\
CEN	&	30.16	\\
TRN	&	59.91	\\
CAN	&	15.80	\\
NEL	&	22.24	\\
OTG	&	40.23	
\end{tabular}
\end{adjustbox}
\caption{Small fixed generation (including geothermal) for each region where this exists.}
\label{tab:fixedgeneration}
\end{table}

\section{Storage trajectories}\label{appendix:storage}

\begin{figure}[htbp]
\centering
\includegraphics[width=0.8\textwidth]{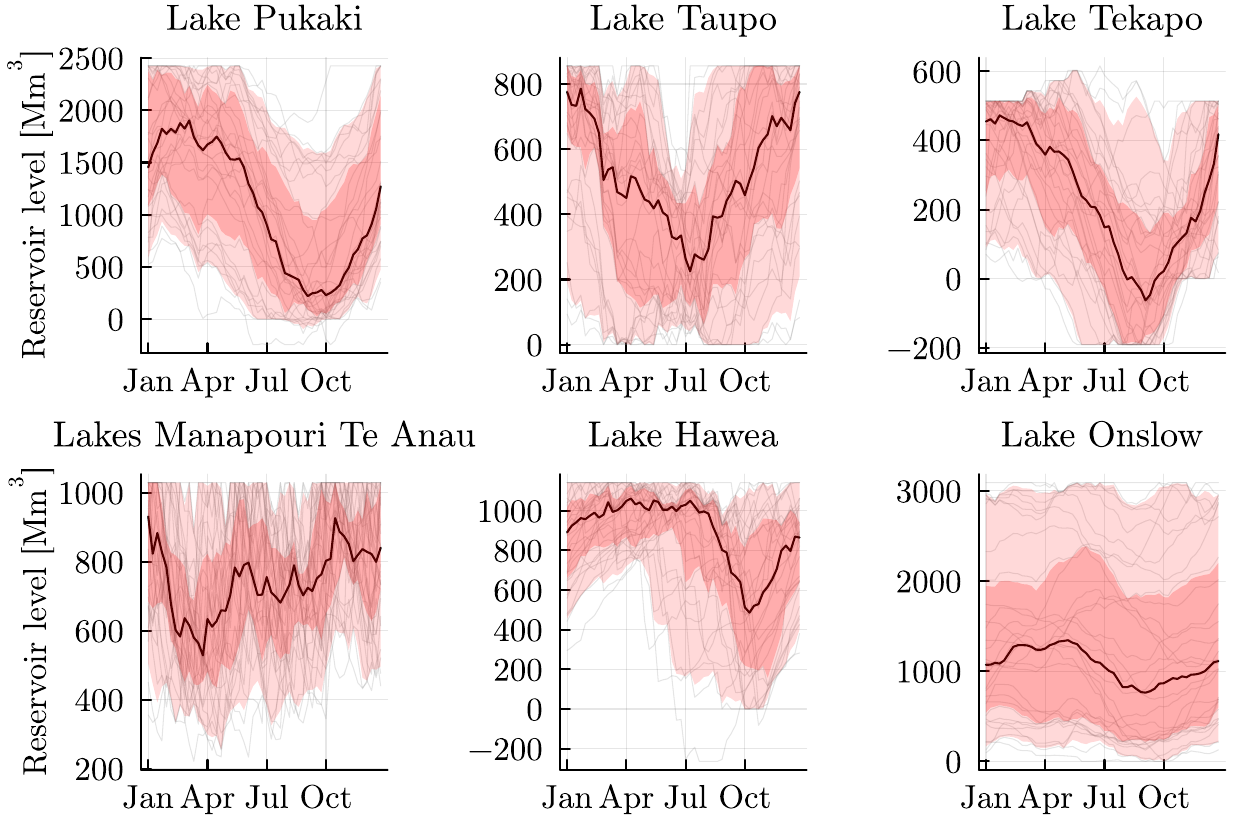}
\caption{Reservoir levels relative to the operating minimum for the historical inflow years 1990--2020 for the \textit{Onslow} case. Light gray lines are the 31 trajectories from each simulated year, the solid line is the median, and the shaded bands are the 10--90 and 25--75 percentiles.  
}
\label{fig:storageonslow1000}
\end{figure}

\begin{figure}[htbp]
\centering
\includegraphics[width=0.8\textwidth]{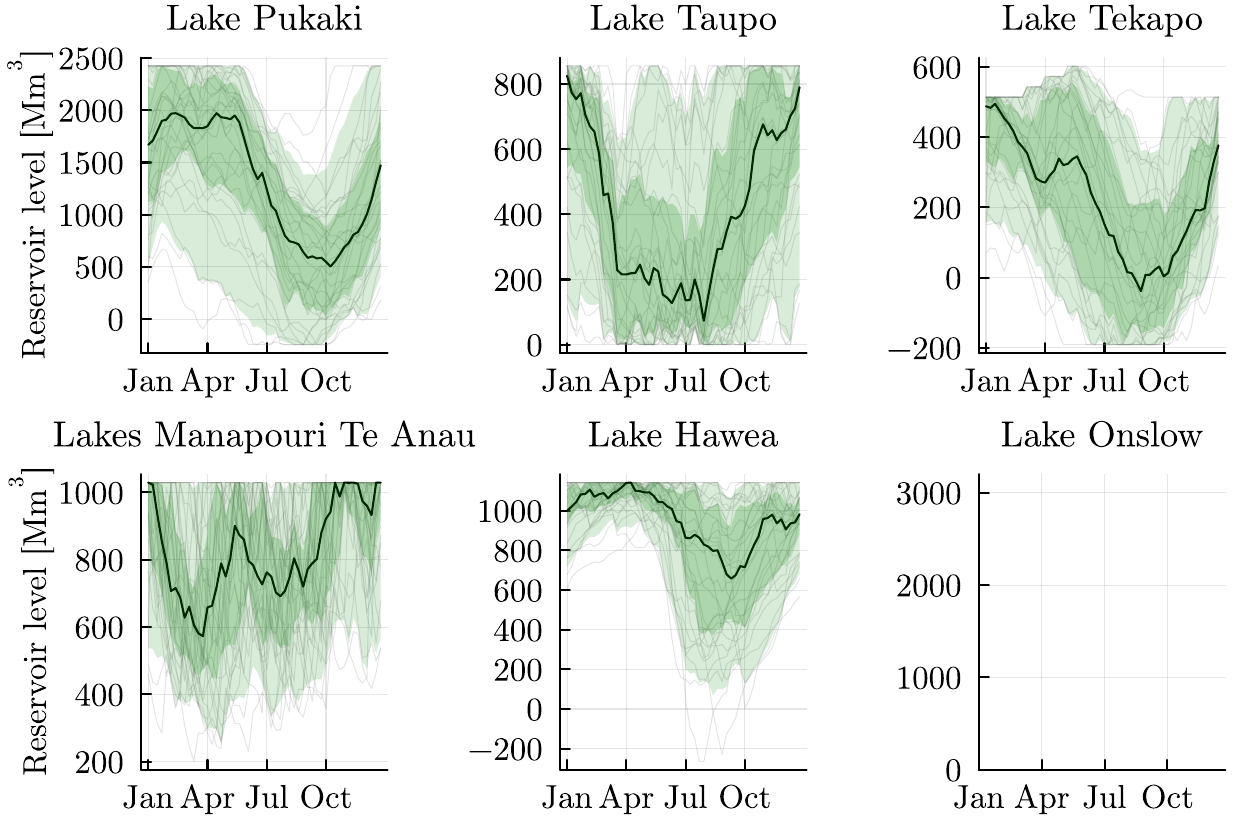} 
\caption{Reservoir levels relative to the operating minimum for the historical inflow years 1990--2020 for the \textit{Peakers} case. Light gray lines are the 31 trajectories from each simulated year, the solid line is the median, and the shaded bands are the 10--90 and 25--75 percentiles.  
}
\label{fig:storagepeakers1000} 
\end{figure}

\begin{figure}[htbp]
\centering
\includegraphics[width=0.8\textwidth]{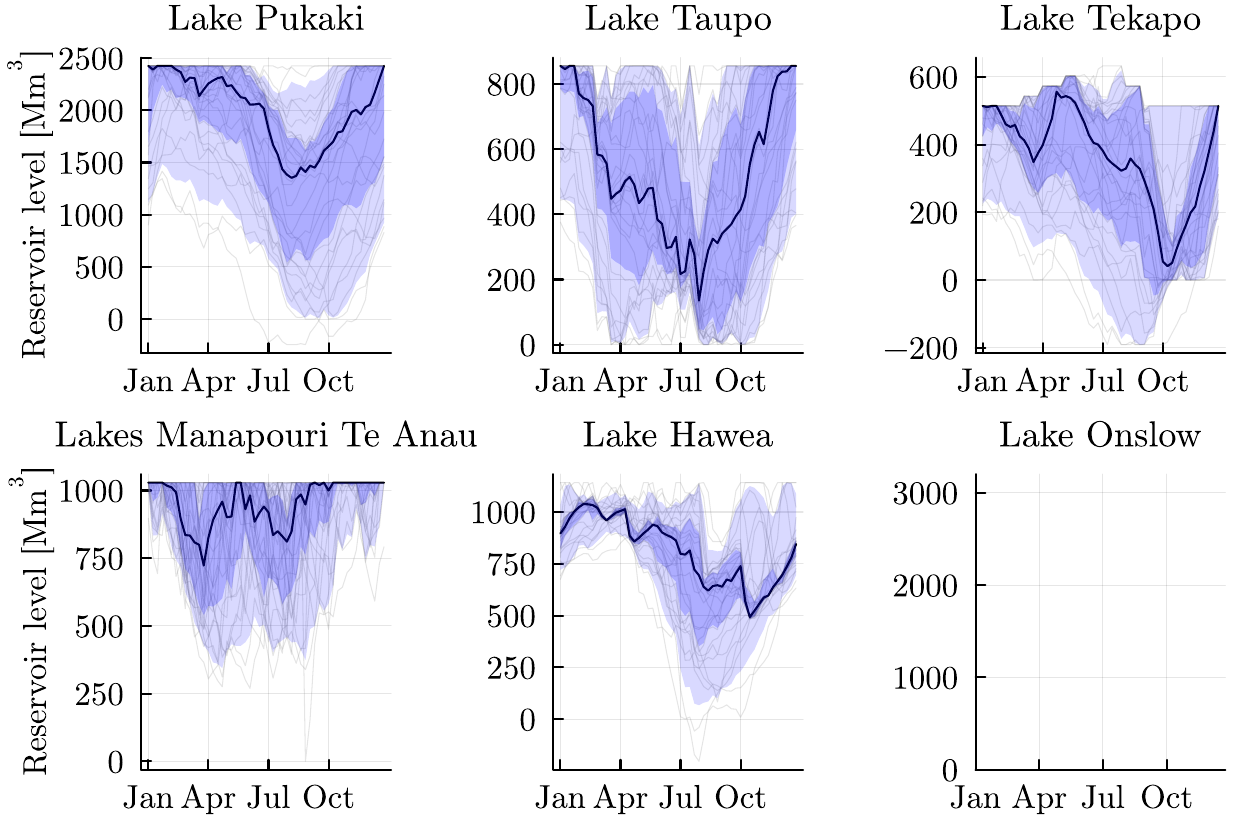} 
\caption{Reservoir levels relative to the operating minimum for the historical inflow years 1990--2020 for the \textit{Wind only} case.
Light gray lines are the 31 trajectories from each simulated year, the solid line is the median, and the shaded bands are the 10--90 and 25--75 percentiles.  
}\label{fig:storagewindonly1000} 
\end{figure}

In Figures \ref{fig:storageonslow1000}, \ref{fig:storagepeakers1000}, and \ref{fig:storagewindonly1000}, we plot the reservoir storage for the six largest reservoirs over the 31 years of historical simulation.

\pagebreak
\section{Marginal water values}\label{appendix:watervalue}

\begin{figure}[htbp]
\centering
\includegraphics[width=0.8\textwidth]{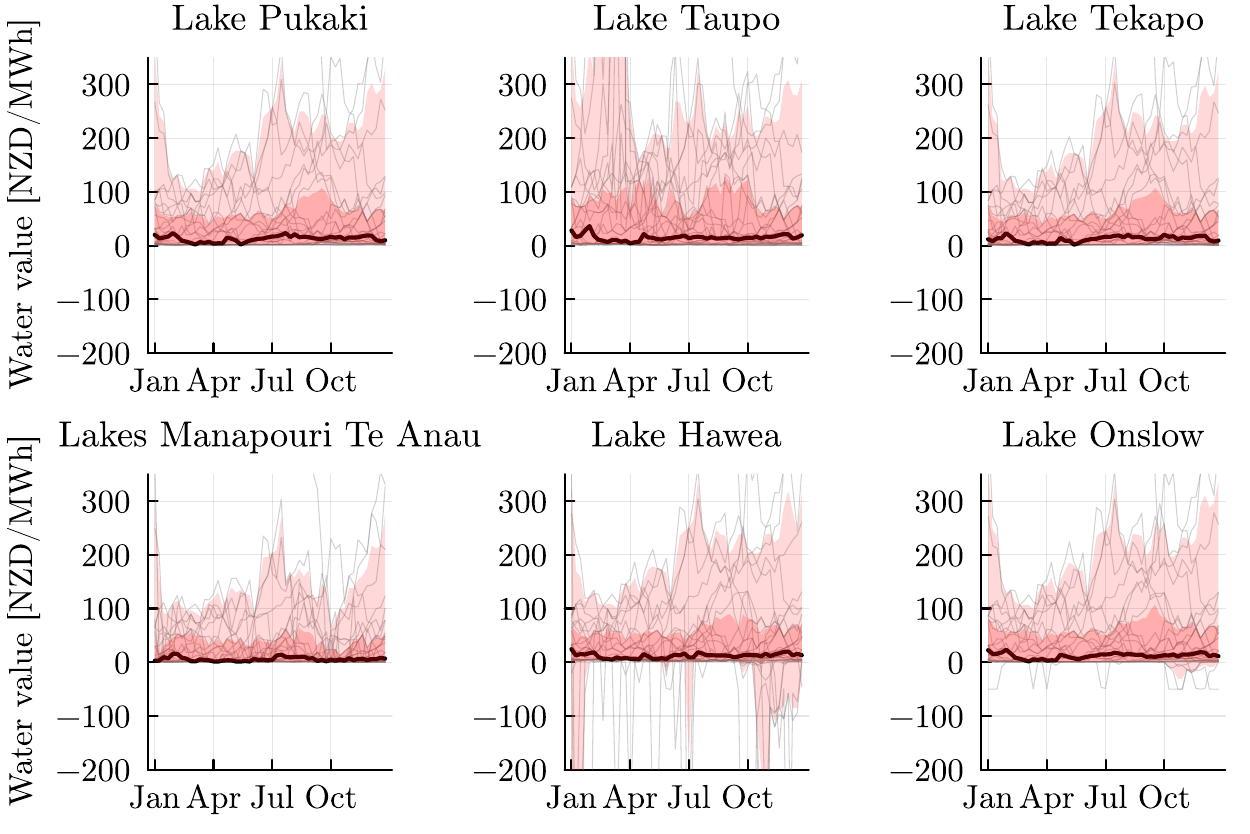}
\caption{Marginal water values for the historical inflow years 1990--2020 for the \textit{Onslow} case. Light gray lines are the 31 trajectories from each simulated year, the solid line is the median, and the shaded bands are the 10--90 and 25--75 percentiles.}
\label{fig:wvonslow1000}
\end{figure}

\begin{figure}[htbb]
\centering
\includegraphics[width=0.8\textwidth]{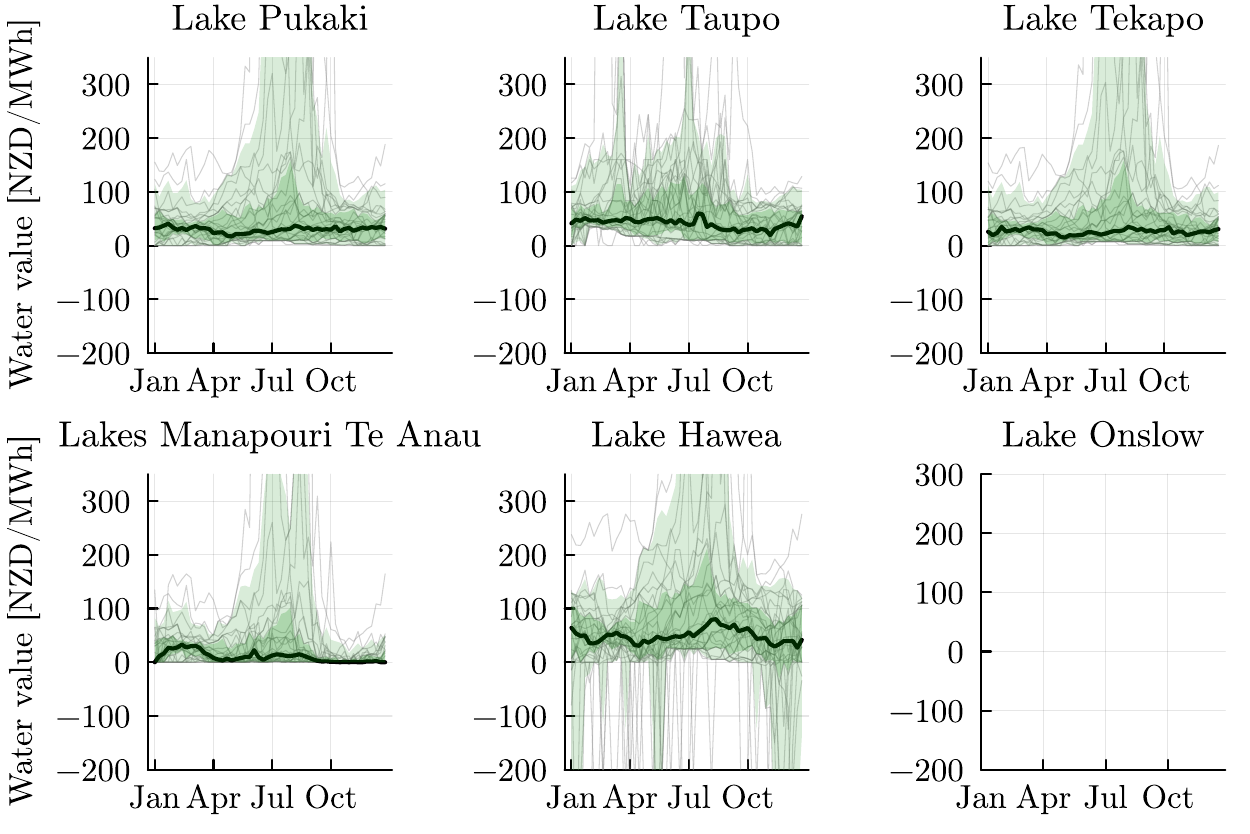} 
\caption{Marginal water values for the historical inflow years 1990--2020 for the \textit{Peakers} case. Light gray lines are the 31 trajectories from each simulated year, the solid line is the median, and the shaded bands are the 10--90 and 25--75 percentiles.}
\label{fig:wvpeakers1000} 
\end{figure}

\begin{figure}[htbp]
\centering
\includegraphics[width=0.8\textwidth]{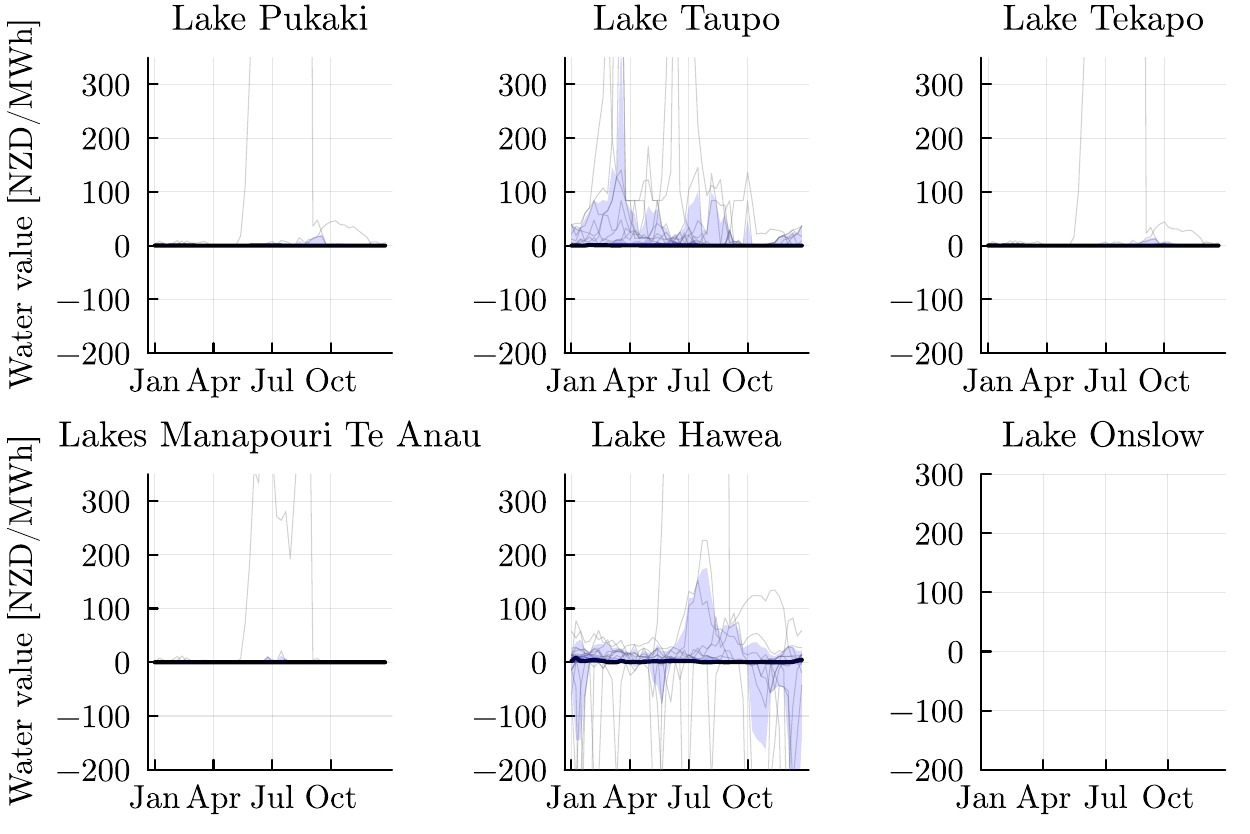} 
\caption{Marginal water values for the historical inflow years 1990--2020 for the \textit{Wind only} case.
Light gray lines are the 31 trajectories from each simulated year, the solid line is the median, and the shaded bands are the 10--90 and 25--75 percentiles.}\label{fig:wvwindonly1000} 
\end{figure}

In Figures \ref{fig:wvonslow1000}, \ref{fig:wvpeakers1000}, and \ref{fig:wvwindonly1000}, we plot the marginal water values for the six largest reservoirs over the 31 years of historical simulation.

\end{document}